\documentclass[11pt,draft]{amsart}

\usepackage{amscd}      
\usepackage{amssymb,amsmath}
\usepackage[all,v2]{xy}
\xyoption{2cell}
\UseAllTwocells
\xyoption{frame}
\CompileMatrices

\usepackage{amsthm}

\newtheorem{prop}{Proposition}[section]
\newtheorem{lem}[prop]{Lemma}

\newtheorem{cor}[prop]{Corollary}
\newtheorem{them}[prop]{Theorem}

\newtheorem{question}[prop]{Question}
\newtheorem{conjecture}[prop]{Conjecture}

\theoremstyle{definition}

\newtheorem{defn}[prop]{Definition}

\newtheorem{rmk}[prop]{Remark}

\newtheorem{numex}[prop]{Example}

\newtheoremstyle{intro}
  {12pt}
  {12pt}
  {\itshape}
  {}
  {\bfseries}
  {.}
  {.5em}
  {}

\theoremstyle{intro}
\newtheorem{introthm}{Theorem}

\newtheorem{intropro}[introthm]{Problem}

\newenvironment{pf}{\begin{trivlist}\item[]{\sc Proof.}}%
            {\nolinebreak $\Box$ \end{trivlist}}
            {\nolinebreak $\Box$ \end{trivlist}}

\newcommand{\noprint}[1]{}

\DeclareMathOperator{\DR}{DR}

\newcommand{\longto}{\longrightarrow}

\newcommand{\B}{{\mathcal B}}

\renewcommand{\tilde}{\widetilde}

\newcommand{\toto}{\rightrightarrows}

\newcommand{\com}{{\scriptscriptstyle\bullet}}
\newcommand{\lcom}{_{\scriptscriptstyle\bullet}}

\newcommand{\upcom}{^{\scriptscriptstyle\bullet}}

\newcommand{\XX}{{\mathfrak X}}

\newcommand{\Gg}{{\mathfrak g}}

\newcommand{\zz}{{\mathbb Z}}

\newcommand{\nn}{{\mathbb N}}

\newcommand{\cc}{{\mathbb C}}
\newcommand{\rr}{{\mathbb R}}

\newcommand{\lL}{\mathcal{L}}

\newcommand{\complex}{{\Bbb C}}

\newcommand{\del}{\partial}

\newcommand{\pr}{\mathop{\rm pr}\nolimits}

\newcommand{\Ad}{\mathop{\rm Ad}\nolimits}

\newcommand{\Der}{\mathop{\rm Der}\nolimits}

\newcommand{\smalcirc}{\mbox{\tiny{$\circ $}}}

\newcommand{\ldiag}[1]%
       {\makebox[0cm]{${\scriptstyle#1}\downarrow\phantom{\scriptstyle#1}$}}
\newcommand{\ldiagup}[1]%
       {\makebox[0cm]{${\scriptstyle#1}\uparrow\phantom{\scriptstyle#1}$}}
\newcommand{\rdiag}[1]%
       {\makebox[0cm]{$\phantom{\scriptstyle#1}\downarrow{\scriptstyle#1}$}}
\newcommand{\sediagr}[1]%
       {\makebox[0cm]{$\phantom{\scriptstyle#1}\searrow{\scriptstyle#1}$}}
\newcommand{\nediagr}[1]%
       {\makebox[0cm]{$\phantom{\scriptstyle#1}\nearrow{\scriptstyle#1}$}}
\newcommand{\rdiagup}[1]%
       {\makebox[0cm]{$\phantom{\scriptstyle#1}\uparrow{\scriptstyle#1}$}}
\newcommand{\swdiag}[1]%
       {\makebox[0cm]{$\phantom{\scriptstyle#1}\swarrow{\scriptstyle#1}$}}
\newcommand{\sediag}[1]%
       {\makebox[0cm]{${\scriptstyle#1}\searrow\phantom{\scriptstyle#1}$}}
\newcommand{\nediag}[1]%
       {\makebox[0cm]{${\scriptstyle#1}\nearrow\phantom{\scriptstyle#1}$}}

\newcommand{\iso}{\stackrel{\sim}{\rightarrow}}

\newcommand{\doublearrowstack}[2]%
 {{{{\scriptstyle#1}\atop{\textstyle\longrightarrow}}\atop{{\textstyle\longright
arrow}\atop{\scriptstyle#2}}}}
\newcommand{\rightleftarrowstack}[2]%
 {{{{\scriptstyle#1}\atop{\textstyle\longrightarrow}}\atop{{\textstyle\longlefta
rrow}\atop{\scriptstyle#2}}}}
\newcommand{\leftrightarrowstack}[2]%
 {{{{\scriptstyle#1}\atop{\textstyle\longleftarrow}}\atop{{\textstyle\longrighta
rrow}\atop{\scriptstyle#2}}}}

\newcommand{\overtoparrow}%
{\makebox[0cm]{\beginpicture
\setcoordinatesystem units <.8cm,.4cm> point at 0 0
\setplotarea x from -3 to 3, y from 0 to 1
\setquadratic
\plot -3 0 0 1 3 0 /
\put{\vector(3,-1){0}}[Bl] at 3 0
\endpicture}}

\newcommand{\underbottomarrow}%
{\makebox[0cm]{\beginpicture
\setcoordinatesystem units <.8cm,.4cm> point at 0 0
\setplotarea x from -3 to 3, y from 0 to 1
\setquadratic
\plot -3 1 0 0 3 1 /
\put{\vector(3,1){0}}[Bl] at 3 1
\endpicture}}

\newcommand{\ses}[5]%
{0\longrightarrow#1\stackrel{#2}{ \longrightarrow}#3\stackrel{#4}{
\longrightarrow}#5\longrightarrow0}

\newcommand{\dt}[6]%
{#1\stackrel{#2}{\longrightarrow}#3 \stackrel{#4}{\longrightarrow}#5
\stackrel{#6}{\longrightarrow} #1[1]}

\newcommand{\cat}[1]%
{(\mbox{\rm #1})}

\newcommand{\gm}{\Gamma}
\newcommand{\tgm}{{\tilde{\Gamma}}}

\newcommand{\tx}{\tilde{x}}

\newcommand{\be }{\begin{eqnarray*}}
\newcommand{\ee }{\end{eqnarray*}}

\newcommand{\edd}{d^\alpha}

\newcommand{\tX}{{\tilde{X}}}

\newcommand{\per}{\backl}
\newcommand{\backl}{\mathbin{\vrule width1.5ex height.4pt\vrule height1.5ex}}

\newcommand{\connections}{connections }

\newcommand{\connection}{connection }
\newcommand{\tdel}{\tilde{\del}}

\newcommand{\frakg}{\Gg}

\newcommand{\cala}{\mathcal{A}}
\newcommand{\calo}{\mathcal{O}}
\newcommand{\bOmega}{\bar{\Omega}}
\newcommand{\ed}{d}
\newcommand{\xxto}[1]{\xrightarrow{#1}}
\newcommand{\ty}{\tilde{y}}
\newcommand{\act}{\star}

\newcommand{\delt}{\widetilde{\del}}
\newcommand{\bg}{B_G}
\newcommand{\OOG}{\OO_G}
\newcommand{\OO}{\Omega} 
\newcommand{\beq}[1]{\begin{equation}\label{#1}}
\newcommand{\eeq}{\end{equation}}

\newcommand{\curv}{\CURV_G(\theta)}

\DeclareMathOperator{\CURV}{curv}
\newcommand{\dg}{d_G}

\newcommand{\etag}{\eta_G}
\newcommand{\eqcls}[1]{\left[#1\right]}
\newcommand{\nablaa}{\nabla}
\newcommand{\PC}{\text{PC}}
\newcommand{\Tr}{{\text{Tr}}}
\newcommand{\nice}{ nice }
\newcommand{\alphab}{\Omega}

\def\gpd{\,\lower1pt\hbox{$\longrightarrow$}\hskip-.24in\raise2pt
             \hbox{$\longrightarrow$}\,}

\newcommand{\gmmt}{\tilde{\Gamma}}
\newcommand{\barom}[1]{(\bOmega\upcom (M, G, L)_{#1}, \edd_{G^{#1}}) }

\title{Periodic Cyclic Homology and Equivariant Gerbes}
\thanks{Research partially supported by the National
 Science Foundation grants DMS 1101827 and DMS 1406668,
and the ANR project ANR-10-BLAN-0111-01-KIND.}

\author{Jean-Louis Tu}
\address{D\'epartement de math\'ematiques, universit\'e de Lorraine}
\email{jean-louis.tu@univ-lorraine.fr}

\author{Ping Xu}
\address{Department of Mathematics, Penn State University}
\email{ping@math.psu.edu}

\begin{document}
\sloppy
\maketitle

\begin{abstract}
This paper is our  first step in   establishing a
de Rham model for equivariant twisted $K$-theory 
using machinery from noncommutative geometry.
Let $G$ be a compact Lie group, $M$ a compact manifold  on
which $G$ acts smoothly. For any $\alpha \in H^3_G (M, \zz)$
we introduce  a notion of localized equivariant  twisted  cohomology
$H^\bullet (\bOmega\upcom (M, G, L)_g,
\edd_{G^g})$,  indexed by $g\in G$.
  We prove that there exists a natural
 family of chain maps, indexed by $g\in G$,
inducing a family of morphisms from 
 the equivariant periodic cyclic homology $HP^G_\com ( C^\infty (M, \alpha ) )$,
where $C^\infty (M, \alpha )$ is a certain smooth algebra
constructed from an equivariant bundle gerbe defined by
$\alpha \in H^3_G (M, \zz)$, to $H^\bullet (\bOmega\upcom (M, G, L)_g,
\edd_{G^g})$. We formulate a conjecture of 
Atiyah-Hirzebruch type theorem for 
equivariant twisted $K$-theory.
\end{abstract}


\section{Introduction}

The well-known Atiyah-Hirzebruch  theorem asserts that for
a smooth manifold $M$, the Chern character  establishes 
an isomorphism: 

$$ K^\bullet(M)\otimes\cc \underset{\cong}{\overset{\text{ch}}{\longto}} H^\bullet_{\DR}(M,\cc)$$ 

Therefore, modulo the torsion, K-theory groups are isomorphic
to the ($\zz_2$-graded) de Rham cohomology groups. 
In the study of  K-theory, it has been a central question how
to establish a Atiyah-Hirzebruch type theorem for other types
of K-theory groups, among which are equivariant K-theory \cite{Segal}.
In 1994, Block-Getzler proved the following remarkable theorem \cite{BG}
extending a result of Baum-Brylinski-MacPherson \cite{BBM} in
the case of $G=S^1$:

\bigskip
Let $G$ be a compact Lie group, $M$ a compact manifold 
on which $G$ acts smoothly. Then
$$K\upcom_{G}(M)\otimes_{R(G)}R^\infty (G)
\underset{\cong}{\overset{\text{}}{\longto}}  
H\upcom (\cala_G\upcom(M), d_{\text{eq}}).$$
\bigskip

Here $R(G)$ is the representation ring
of $G$,  and $R^\infty (G)$ is the ring of
smooth functions on $G$ invariant under the conjugation.
Then $R^\infty (G)$ is an algebra over $R(G)$,
 since $R(G)$ maps to
$R^\infty (G)$ by the character map.
And $H\upcom (\cala_G\upcom(M), d_{\text{eq}})$ is
the cohomology of global equivariant  differential 
forms on $M$, also called the {\em delocalized equivariant} 
cohomology. Roughly speaking,
 $H\upcom (\cala_G\upcom(M), d_{\text{eq}})$
can be considered as  the cohomology of the  inertia  stack
 $\Lambda [M/G]$ (while ordinary  equivariant cohomology
is the cohomology of the quotient stack $[M/G]$).
 In a certain sense, delocalized equivariant
cohomology  $H\upcom (\cala_G\upcom(M), d_{\text{eq}})$ is a
de Rham description of the equivariant  $K$-theory.
 
In late 1980's,  Block \cite{Block} and Brylinski \cite{Bry1, Bry2}
 independently proved the following theorem:

\bigskip
Let $G$ be a compact Lie group and  $A$  a topological  $G$-algebra.
 Then the equivariant Chern character \eqref{eq:eqchern} induces an isomorphism
\begin{equation}
HP_\bullet ^G(A) \underset{\cong}{\overset{\text{}}{\longto}} K_\bullet^G(A)\otimes_{R(G)}R^\infty (G).
\end{equation}
\bigskip

 By using the above theorem,  Block-Getzler  \cite{BG}
reduces the problem of establishing  Atiyah-Hirzebruch type theorem for 
equivariant K-theory  to that of computing
the  equivariant periodic cyclic homology $HP_\bullet^G (C^\infty (M))$,
and thus can apply  the machinery of noncommutative geometry.
When $G=\{*\}$, $HP_\bullet^G (C^\infty (M))$ is isomorphic to
the ($\zz_2$-graded) de Rham cohomology $H_{DR}^\bullet (M)$
 according to a theorem of Connes
\cite{Con85, Con:book}, and therefore Block-Getzler theorem
reduces to the classical Atiyah-Hirzebruch  theorem.

Motivated by string theory, there has been a great deal of
interest in the study of twisted $K$-theory \cite{BM, Witten}.
 It is thus natural to ask how to extend Atiyah-Hirzebruch  theorem
to twisted  $K$-theory.  The    manifolds case
was completely solved in \cite{Mat-Ste04}, while  orbifolds case
was done \cite{TX:06}. 
 This paper is
our first step  in establishing Atiyah-Hirzebruch type  theorem
for twisted equivariant $K$-theory.

 In literatures, there 
exist different equivalent approaches to twisted equivariant $K$-theory,
eg. \cite{Ati-Seg}. In \cite{TXL:04}, with Laurent-Gengoux,
we introduced twisted equivariant $K$-theory  based on
the idea from noncommutative geometry. It can be described
roughly as follows. For  a compact manifold $M$
equipped with an action of  a compact Lie group $G$,
and   any   $\alpha \in H^3_G(M, \zz )$, one can always
construct an $S^1$-central extension of Lie groupoids
$\tilde{X}_1\stackrel{\pi}{\to}X_1 \toto X_0$ 
representing $\alpha \in H^3_G (M, \zz)$. Such an $S^1$-central extension
is unique up to Morita equivalence. From the $S^1$-central extension
 of Lie groupoids $\tilde{X}_1\stackrel{\pi}{\to}X_1 \toto X_0$,
one constructs a convolution algebra $C_c(X_1, L)$. Then the twisted 
equivariant $K$-theory groups can be defined as the
$K$-theory groups of this  algebra (or its corresponding
reduced $C^*$-algebra), i.e.,
 $K\upcom_{G, \alpha}(M):=K_\bullet (C_c(X_1, L))$.

In order to apply Block-Brylinski theorem \cite{Block, Bry1, Bry2}
(Theorem \ref{thm:Block-Bry})
in our situation, first of all,  we prove the following

\begin{introthm}\label{AAA}
For  any integer class $\alpha \in H^3_G(M, \zz)$,
there always  exists a $G$-equivariant bundle gerbe \cite{Stienon}
 $\tilde{H}_1\stackrel{\phi}{\to} H_1\toto H_0$
over  $M$  with an equivariant connection and
an equivariant curving, whose
 equivariant $3$-curvature
represents $\alpha$ in the Cartan model $(\Omega^\bullet_G(M), d+ \iota)$.
\end{introthm}

As a consequence,  equivalently  one can    define the twisted
equivariant $K$-theory  groups $K\upcom_{G, \alpha}(M)$
as $K_\bullet^G(C_c^\infty(H,L))$, where $C_c^\infty(H,L)$
is the convolution algebra corresponding to the $G$-equivariant bundle gerbe
 $\tilde{H}_1\stackrel{\phi}{\to} H_1\toto H_0$.
Note that when $G=\{*\}$ and $\alpha \in H^3_G(M, \zz )$ being trivial,
$\tilde{H}_1\stackrel{\phi}{\to} H_1\toto H_0$ is Morita equivalent
to $C^\infty (M)$.

Therefore, we are led to  the following 

\begin{intropro}
\label{problem:A}
Compute the  equivariant periodic cyclic homology
 $HP_\bullet^G (C_c^\infty(H,L))$
in terms of geometric data to obtain de Rham type cohomology groups.
\end{intropro}

Toward this direction, we first prove the following:

\begin{introthm}
Any Lie groupoid $S^1$-central extension  representing 
$\alpha \in H^3_G(M, \zz )$ canonically induces a family of  $G$-equivariant
flat  $S^1$-bundles $\coprod_{g\in G}
(P^g\to M^g)$ indexed by $g\in G$, where 
$M^g=\{x\in M|x\cdot g=x\}$ is the fixed point set under the
diffeomorphism $x\to x\cdot g$, and for any $h\in G$,  the action by $h$,
  is an isomorphism of  the flat bundles $(P^g\to M^g)  \mapsto
 (P^{h^{-1}gh}\to M^{h^{-1}gh})$.

Two such  $S^1$-central extensions induce $G$-equivariantly isomorphic
families of  flat  $S^1$-bundles.
\end{introthm}

Let $L=\coprod_{g\in G} L^g$,
where $L^g=P^g\times_{S^1}\complex$, $\forall g\in G$,
are their associated  $G$-equivariant flat complex line bundles.

Choose a $G$-equivariant closed $3$-form $\eta_G\in Z^3_G (M)$
such that $[\eta_G]=\alpha$.
Following Block-Getzler \cite{BG}, we
consider the  localized twisted equivariant cohomology as follows.
Denote by $\bOmega\upcom (M, G, L)_g$ the space $\bOmega_{G^g} (M^g, L^g )$ of germs at zero of $G^g$-equivariant
 smooth maps from $\frakg^g$ to $\Omega^\com (M^g , L^g)$.
By $\edd_{G^g}$ we  denote the twisted  equivariant differential operator
$\nablaa_g+ \iota-2\pi i \eta_G$ on $\bOmega_{G^g} (M^g , L^g)$,
where $\nablaa_g: \bOmega_{G^g} (M^g , L^g)\to \bOmega_{G^g} (M^g , L^g)$
is  the covariant differential  induced by  the flat connection
on the complex line bundle $L^g\to M^g$, and $\eta_G$ acts
on $\bOmega_{G^g} (M^g , L^g)$ by  taking the wedge product    with
$i_g^* \eta_G$. Here $i_g^*: \Omega_{G} (M)\to \Omega_{G^g} (M^g )$
is the restriction map. It is simple to see that $(\edd_{G^g})^2=0$,
and  $\{ \edd_{G^g}|g\in G\}$ are compatible with the $G$-action.
The family of cohomology groups
are denoted by $H^\bullet (\bOmega\upcom (M, G, L)_g, \edd_{G^g})$, 
and are called {\em localized twisted equivariant cohomology}.

Our main result in the paper is the following:

\begin{introthm}
There exists a family of $G$-equivariant chain maps, indexed by   $g\in G$:
$$\tau_g:(\PC_\bullet^G(C_c^\infty(H,L)),b+\B)
\to  \barom{g}
$$
By  $G$-equivariant chain maps, we mean that
the following diagram
of chain maps commutes:
$$\xymatrix{
(\PC_\bullet^G(C_c^\infty(H,L)),b+\B) \ar[r]^{\tau_{g}}
\ar[dr]_{\tau_{h^{-1}gh}} & ( \barom{g})
\ar[d]^{\phi_{h}}\\
& ( \barom{h^{-1}gh}).
}$$

Therefore, there is a family of morphisms on the level of cohomology:
\begin{equation}
\tau_g: HP^G_\bullet (C_c^\infty(H,L)) \to   H^\bullet (  \barom{g}).
\end{equation}
\end{introthm}

In order to completely  solve Problem \ref{problem:A},
 following Block-Getzler \cite{BG}, we  propose the following

\begin{intropro}
\label{problem:B}
Introduce global twisted equivariant
differential forms by modifying   the notion of global equivariant
differential forms a la Block-Getzler \cite{BG} to
 define  delocalized twisted equivariant cohomology $H_{G, delocalized, \alpha}\upcom (M)$,
and establish the isomorphism
 $$ HP_\bullet^G (C_c^\infty(H,L)) \underset{\cong}{\overset{\text{}}{\longto}}
 H_{G, delocalized, \alpha}\upcom (M). $$ 
\end{intropro}

We will  devote Section \ref{sec:discussion} to 
 discussions on this issue.

\subsection*{Acknowledgments}
We would like to thank several institutions for their hospitality
 while work on this project was being done:
Penn State University (Tu),
IHES (Xu),  and Universit\'e Paris Diderot (Xu).
We also wish to thank many people for useful discussions and comments,
 including  Mathieu Sti\'enon and Camille Laurent-Gengoux.

\section{Localized equivariant twisted  cohomology}

\subsection{$S^1$-gerbes over $[M/G]$}

Let  $\tX_1 \to X_1\toto X_0$ be an $S^1$-central extension
of Lie groupoids.  By abuse of notations,
  we  denote, by  $d$,  both the de Rham differentials
$\Omega^\bullet (X_\com )\to \Omega^{\bullet+1} (X_\com )$,
and $\Omega^\bullet (\tX_\com )\to \Omega^{\bullet+1} (\tX_\com )$.
And, by $\del$, we denote the simplicial differential  
$\del: \Omega^\bullet (X_\com ) \to \Omega^\bullet (X_{\bullet +1})$
 for the groupoid $X_1\toto X_0$,
while,  by $\tdel$, we denote the  simplicial differential
$\tdel :
\Omega^\bullet (\tX_\com ) \to \Omega^\bullet (\tX_{\bullet +1})$
 for the groupoid $\tX_1 \toto X_0$. See Section 2.1 \cite{Xu:JDG04}
for  details   on de Rham cohomology of Lie groupoids.
  Recall  the following \begin{defn}[ \cite{BX}]

(i) A connection form $\theta\in \Omega^1(\tX_1)$
for the  $S^1$-bundle $\tX_1\to X_1$, such that
$\tdel\theta=0$,
 is a {\em \connection} on the  $S^1$-central extension $S^1\to
\tX_1 \to X_1\toto X_0$;

(ii) Given $\theta$, a 2-form $B\in\Omega^2(X_0 )$, such that $d\theta=\tdel B$
is a {\em curving};

(iii) and given $(\theta,B)$, the 3-form
$\Omega=dB\in H^0(X\lcom,\Omega^3)\subset \Omega^3(X_0 )$ is called
the {\em $3$-curvature}.

\end{defn}

We have the following

\begin{lem}
Given an $S^1$-central extension  $\tX_1\stackrel{\pi}{\to} X_1\toto X_0$,
\begin{enumerate}
\item the obstruction group to the existence of \connections
is $H^2(X\lcom , \Omega^1)$;
\item the obstruction group to the existence of curvings
is $H^1(X\lcom , \Omega^2)$.
\end{enumerate}
\end{lem}
\begin{pf}
(1). 
Choose any connection one-form $\theta\in \Omega^1 (\tX_1 )$ of
the $S^1$-bundle $\pi: \tX_1\to X_1$. It is simple to see that 
$\tdel \theta=\pi^* \eta$, where $\eta \in \Omega^1 (X_2)$.
Here by abuse of notations, we use the same symbol $\pi$
to denote the induced projection $\tX_2\to X_2$. 
Since $\pi^* \del \eta=\tdel \pi^*  \eta=\tdel^2 \theta=0$,
it  thus follows that  $\del \eta=0$. If $\theta'$ is another connection
one-form, then $\theta'$ differs from $\theta$ by a 
one-form $A\in \Omega^1(X_1)$.
Thus $\eta'=\eta+\del A$. It follows that the class
 $[\eta]\in H^2 (\Omega^1(X\lcom ), \del)$
is independent of  the choice of $\theta$.
Note that $H^\bullet (\Omega^1(X\lcom ),\del)\cong H^\bullet
 (X\lcom,\Omega^1)$ since $\Omega^1$ is a soft sheaf.

Assume that  $[\eta]\in H^2(X\lcom , \Omega^1)$ vanishes.
We may write $\eta =\del \alpha $ for some one-form $\alpha \in
\Omega^1( X_1 )$. It is simple to see that $\theta'=\theta-\pi^*
\alpha$ is  indeed a \connection for the $S^1$-extension.

(2). Assume that $\theta\in \Omega^1 (\tX_1 )$ is
a  connection. Let $\omega\in \Omega^2 (X_1)$
be its curvature, i.e. $d\theta =\pi^* \omega$.
Since $\pi^*\del\omega=\tdel \pi^* \omega =\tdel d\theta=d \tdel \theta=0$,
 we have $\del\omega=0$. 
Hence  $[\omega]\in H^1 (\Omega^2(X\lcom ),\del)\cong
 H^1(X\lcom , \Omega^2)$.  
Then $[\omega]$ vanishes
 if and only if there exists $B\in\Omega^2(X_0 )$ such that $\omega 
=\del B$, i.e. $d\theta=\tdel B$. Hence $[\omega]=0$ if and only
if there exists a curving.
\end{pf}

\begin{rmk}
Note that  $H^2(X\lcom , \Omega^1)$ and 
$H^1(X\lcom , \Omega^2)$ are isomorphic to
$H^2 (\XX, \Omega^1)$ and $H^1 (\XX, \Omega^2)$, 
 respectively,
where $\XX$ is the  differentiable   stack corresponding to the
groupoid  $X_1\toto X_0$ \cite{BX}.
Therefore, these cohomology groups are Morita invariant.
\end{rmk} 

The following theorem is due to Abad-Crainic [\cite{AC} Corollary 4.2].

\begin{them}
If  $X_1\toto X_0$ is a proper Lie groupoid,
then  
$$H^p (X_\com, \Omega^q)=0, \ \ \ \mbox{if } p>q.$$ 
\end{them}

In particular, we have $H^2 (X_\com , \Omega^1)=0$. Thus we have 
the following

\begin{prop}
\label{pro:connection}
If  $ \tX_1 \to X_1\toto X_0$ is  an $S^1$-central extension
of a proper Lie groupoid $X_1\toto X_0$,
then connections  always exist.

In particular, if
$G$ is a compact Lie group acting on a manifold $M$, and
$\tilde{X}_1\stackrel{\pi}{\to}X_1 \toto X_0$ is  an
$S^1$-central extension representing any class
$\alpha \in H^3_G (M, \zz)$, then this extension admits   a connection.
\end{prop}

Note that  $H^1 ((M\rtimes G)\lcom, \Omega^2)$ may not 
necessarily  vanish in general, so curvings does not always exist.

\subsection{Equivariant bundle gerbes}

We now recall some basic notions regarding equivariant cohomology
in order to fix the notations.

Let $M$ be a smooth 
manifold with a smooth right action of a compact Lie group $G$:
$(x, g)\in M\times G\to x\cdot g$.
There is an induced action of the group $G$  on
the space $\Omega^\bullet (M)$ of differential forms on $M$ by
$g\cdot \omega=R_g^* \omega$, where $R_g: M\to M$ is the
operation of the action by $g\in G$.
If $\omega: \frakg \to \Omega^\bullet (M)$
is a  map from $\frakg$  to $\Omega^\bullet (M)$, the group $G$
 acts on $\omega$ by the formula
\begin{equation}
(g\cdot \omega)(X)=g\cdot  (\omega (\Ad_{g^{-l}} X)), \ \ \ \forall X\in \frakg.
\end{equation}

By an  equivariant differential form on $M$, we mean
a $G$-equivariant polynomial function $\omega: \frakg \to \Omega^\bullet (M)$.
When $\omega$ is a homogeneous polynomial, the degree of $\omega$
 is defined to be the sum of $2\times$the degree
of the polynomial and  the degree of the differential form $\omega (X), \ X\in
\frakg$.  We denote by  $\Omega^k_G(M)$ the space of
local equivariant differential forms of degree $k$.
Define an equivariant differential $d_G=d+  \iota:
\Omega^\bullet_G(M)\to \Omega^{\bullet+1}_G(M)$, where
$$(d\omega )(X)=d(\omega (X)), \ \ \ (\iota \omega)(X)=\iota_{\hat{X}} \omega (X). $$
Here $\hat{X}$ denotes
  the infinitesimal vector field on $M$ generated by
the Lie algebra element $X\in \frakg$.  
Note that  $\Omega^\bullet_G(M)$ can also be 
identified with  the space of invariant polynomials
$(S\frakg^* \otimes \Omega^\com (M))^G$.
 The cochain complex  $(\Omega^\bullet_G(M), d_G)$ is 
called the  Cartan model of the
  equivariant cohomology group $H_G^\bullet (M)$ \cite{BGV}.


Following \cite{Stienon},
  an $S^1$-central extension $\tilde{H}_1\xxto{p}H_1\toto H_0$
 is said to be
$G$-equivariant if both $\tilde{H}_1\toto H_0$ and
 $H_1\toto H_0$ are $G$-groupoids, the groupoid morphism
 $p:\tilde{H}_1\to H_1$ is $G$-equivariant,
 and $\tilde{H}_1\xxto{p} H_1$ is a $G$-equivariant principal $S^1$-bundle, i.e.
 if the following relations:
\begin{gather*}
(\tx\cdot\ty)\act g=(\tx\act g)\cdot(\ty\act g) \\
p(\tx\act g)=p(\tx)\act g  \\
(\lambda\tx)\act g=\lambda(\tx\act g)
\end{gather*}
are satisfied for all $g\in G$, all composable pairs $(\tx,\ty)$ in $\tilde{H}_2$ and all $\lambda\in S^1$.

Now assume that  $H_0 \xxto{\pi}M$ is  a $G$-equivariant surjective submersion.
Consider the pair groupoid $ H_1\toto H_0$, where $H_1=H_0\times_M H_0$,
the source and target maps are 
$t(x,y)=x$ and $s(x,y)=y$, and the multiplication $(x,y)\cdot(y,z)=(x,z)$. 
Then $ H_1\toto H_0$ is a 
$G$-groupoid, which is Morita equivalent to the $G$-manifold
$M\toto M$. A {\em $G$-equivariant bundle gerbe}
 \cite{Meinrenken:gerbe, Stienon} is a $G$-equivariant
$S^1$-central extension of Lie groupoids
$\tilde{H}_1\xxto{p}H_1\toto H_0$.

The following notion is due to Stienon \cite{Stienon}.

\begin{defn} [\cite{Stienon}]
\label{curving}
\begin{enumerate}
\item An equivariant connection is a $G$-invariant 1-form $\theta\in\OO^1(\tilde{H}_1 )^G$ 
such that $\theta$ is a connection 1-form for the principal $S^1$-bundle
 $\tilde{H}_1\xxto{p}H_1 $ and satisfies
\[ \delt\theta=0 .\]
\item Given an equivariant connection $\theta$, an equivariant curving 
is a degree-2 element $\bg\in\OOG^2(H_0 )$ such that
\beq{2} \curv=\del\bg ,\eeq
where $\curv$ denotes the equivariant curvature of the $S^1$-principal 
bundle $\tilde{H}_1\xxto{p}H_1$, i.e. 
the element $\curv\in\OOG^2(H_1)$ characterized by the 
relation
\beq{1} \dg\theta=p^*\curv .\eeq
\item Given an equivariant connection and an equivariant
curving $(\theta,\bg)$, the corresponding equivariant 3-curvature 
is the equivariant 3-form $\etag\in\OOG^3(M)$ such that
\beq{3} \pi^*\etag=\dg\bg .\eeq
\end{enumerate}
Here the coboundary operators
$\del: \Omega^\bullet (H_\com)\to \Omega^\bullet (H_{\com +1})$
and $\delt: \Omega^\bullet (\tilde{H}_\com)\to \Omega^\bullet 
(\tilde{H}_{\com +1})$ are the  simplicial differentials
of  the groupoids $H_1\toto H_0$ and $\tilde{H}_1\toto H_0$,   respectively.
\end{defn}

The following result can be found in \cite{Stienon}.

\begin{prop}
Let $\tilde{H}_1\xxto{\phi }H_1\toto H_0$ be a $G$-equivariant bundle gerbe over a $G$-manifold $M$.
\begin{enumerate}
\item Equivariant connections and equivariant
 curvings $(\theta,\bg)$ always exist.
\item The class $\eqcls{\etag}\in H_G^3(M)$ defined by
the equivariant 3-curvature is independent of the choice of $\theta$ and $\bg$.
\end{enumerate}
\end{prop}

The degree $3$  equivariant cohomology
  class $\eqcls{\etag}\in H_G^3(M)$ is called the
{\em  equivariant Dixmier-Douady class}
 of the underlying equivariant bundle gerbe.

\subsection{Kostant-Weil theorem for equivariant   bundle gerbes}

The main result of this section is the following Kostant-Weil type of 
quantization theorem for equivariant   bundle gerbes.

\begin{them}
\label{thm:main-bundlegerbe}
For  any integer class $\alpha \in H^3_G(M, \zz)$,
there always  exists a $G$-equivariant bundle gerbe
over  $M$  with an equivariant connection and 
an equivariant curving, whose
 equivariant $3$-curvature
represents $\alpha$ in the Cartan model. 
Moreover, one can choose a $G$-basic 
connection as  an equivariant connection on the $G$-equivariant bundle gerbe.
\end{them}

The rest of this subsection is devoted
to  the proof of  Theorem \ref{thm:main-bundlegerbe}.

For any   $\alpha \in H^3_G(M, \zz)  \cong H^3 ((M\rtimes G)\lcom, \zz )$,
let $\tilde{\gm}\stackrel{\pi}{\to}\gm \toto M'$, where
 $p: M'\to M$ is an immersion, be
 an $S^1$-central extension  of Lie groupoids representing
$\alpha$.
Let $H_0=M'\times G$. Assume that the map
$$\pi : H_0 \to M, \ \   \pi (x,g)=p(x)g$$
is a surjective submersion. Note  that this assumption
always holds if $M' =\coprod U_i$ is an open cover of $M$.
  It is clear that $H_0$ admits  a  right $G$-action: 
 $ (x, h)\cdot g=(x, hg), \  \forall (x, h)\in M'\times G$
and $g\in G$. 
Then the map  $\pi : H_0 \to M$ is clearly $G$-equivariant.
Let $H_1=H_0 \times_M H_0$. Then $H_1\toto H_0$ is a $G$-groupoid, which
is Morita equivalent to $M\toto M$.
 Consider the groupoid morphism
 $$\nu: H_1\to \gm,\ \   \nu((x,g),(y,h))=(x,gh^{-1}, y).$$
Let $\tilde{H}_1\to H_1$ be the pullback $S^1$-bundle of
$\tgm\to \gm$ via the map $\nu: H_1\to \gm$:
$$\tilde{H}_1=H_1\times_{\gm }\tgm=\{((x,g),(y,h), \tilde{r})|
(x,g),(y,h)\in H_0,
 \tilde{r}\in \tgm, p(x)g=p(y)h, \ (x, gh^{-1}, y)=p (\tilde{r})\}. $$
Since $\nu$ is constant along the $G$-orbits on $H_1$, the $G$-action on   
$H_1$ naturally  lifts to $\tilde{H}_1$, i.e.
$$((x,g),(y,h), \tilde{r})\cdot k=((x,gk),(y,hk), \tilde{r}), \ \ \forall 
k\in G, ((x,g),(y,h), \tilde{r})\in \tilde{H}_1.$$
Thus  $\tilde{H}_1\toto H_0$ is a $G$-equivariant
groupoid. Let $\phi :\tilde{H}_1\to H_1$ be the projection.

\begin{lem}
\label{pro:H}
\begin{enumerate}
\item  $\tilde{H}_1\stackrel{\phi}{\to} H_1\toto H_0$ is a $G$-equivariant
bundle gerbe over $M$;
\item the following diagram 
\begin{equation}
\label{eq:morphism}
 \xymatrix{
 {\tilde{H}_1}\ar[d]^{} \ar[r]^{\tilde{q}}
& {\tilde{\gm}}
\ar[d]^{} \\
H_1  \ar@<0.5ex>[d] \ar@<-0.5ex>[d] \ar[r]^{q} & \gm
\ar@<0.5ex>[d] \ar@<-0.5ex>[d] \\
H_0 \ar[r]_{q_0 } & M'},
\end{equation}
where
 $\tilde{q}: \tilde{H}_1\to \tilde{\gm}$, $
((x,g),(y,h), \tilde{r})\mapsto \tilde{r}$ is the projection,
$q: H_1\to \gm$ is the map $\nu$,
and $q_0 : H_0  (\cong M' \times G)\to M'$ is the projection
$pr_1$,
defines a morphism of $S^1$-central extensions;
\item the $S^1$-central extensions $\tilde{H}_1\rtimes G
\stackrel{\phi }{\to} H_1\rtimes G \toto H_0$
and  $\tilde{\gm}\stackrel{\pi}{\to}\gm \toto M'$
are Morita equivalent. Indeed,
\begin{equation}
\label{eq:morita}
 \xymatrix{
 {\tilde{H}_1}\rtimes G\ar[d]^{} \ar[r]^{\tilde{q} \smalcirc  pr_1} 
& {\tilde{\gm}}
\ar[d]^{} \\
H_1\rtimes G  \ar@<0.5ex>[d] \ar@<-0.5ex>[d] \ar[r]^{q\smalcirc  pr_1} & \gm
\ar@<0.5ex>[d] \ar@<-0.5ex>[d] \\
H_0 \ar[r]_{q_0 } & M' }
\end{equation}
is a Morita morphism of $S^1$-central extensions \cite{TXL:04}.
\end{enumerate}
\end{lem}
\begin{pf}
	This can  be verified directly, which
is left to  the reader.
\end{pf}

As an immediate consequence, we have the following

\begin{cor}
\label{pro:1.9}
For  any integer class $\alpha \in H^3_G(M, \zz)$, there exists a 
$G$-equivariant bundle gerbe  $\tilde{H}_1\to H_1\toto H_0$ 
over $M$ such that the Dixmier-Douady
class of $\tilde{H}_1\rtimes G
\stackrel{\phi }{\to} H_1\rtimes G \toto H_0$ is equal to
$\alpha$ under the natural 
 isomorphism 
$$H^3((H\rtimes G)_\com, \zz)\cong H^3((M\rtimes G)_\com, \zz)
 \cong H^3_G(M, \zz).$$
\end{cor}
\begin{pf}
For any   $\alpha \in H^3_G(M, \zz)  \cong H^3 ((M\rtimes G)\lcom, \zz )$,
according to Proposition 4.16 in \cite{BX},   one can always
 represent $\alpha$
by an $S^1$-central extension \cite{BX, Tu:2005}
$\tilde{\gm}\stackrel{\pi}{\to}\gm \toto M'$, where
 $p: M'\to M$ is \'etale and surjective,
 and 
$\gm=\{(x,g,y)\in M'\times G\times M'\vert\;
p(x)g=p(y)\}$ is the pull-back groupoid of $M\rtimes G\toto M$ via
$p$. Then
\begin{equation}
\xymatrix{
\gm \ar[r]^p\ar[d]\ar@<-.6ex>[d] &
 M\rtimes G\ar@<-.6ex>[d]\ar[d] \\
M' \ar[r]^{p} &  M
}
\end{equation}
is a Morita morphism, where,
by abuse of notations, $p: \gm \to M\rtimes G$ is  given by
$p(x,g,y) =(p(x), g)$. The groupoid structure on
$\gm \toto M'$ is given by $s (x,g,y)=y, \ t (x,g,y) =x$,
$(x, g, y)(y, h, z)=(x, gh, z)$ and $(x, g, y)^{-1}=(y, g^{-1}, x)$.
Thus we are in the situation described  at the beginning of this section.
In particular, the Dixmier-Douady class of the $S^1$-central extension
$\tilde{H}_1\rtimes G
\stackrel{ }{\to} H_1\rtimes G \toto H_0$
is equal to $\alpha$ under the isomorphism
$H^3((H \rtimes G)_\com, \zz)\cong H^3((M\rtimes G)_\com, \zz)
=H^3_G(M, \zz)$. This concludes the proof of
the corollary.
\end{pf}

{\bf Proof of Theorem \ref{thm:main-bundlegerbe}}
It remains to prove that the Dixmier-Douady class, in the Cartan model,
 defined by the equivariant $3$-curvature
of the $G$-equivariant bundle gerbe $\tilde{H}_1\to H_1\toto H_0$
in Corollary  \ref{pro:1.9} is equal to $\alpha$. This
essentially follows from  a theorem of Stienon (Theorem 4.6 in
 \cite{Stienon}).
More specifically, if $\eta_G\in Z^3_G(M)$ is the
equivariant $3$-curvature of the $G$-equivariant bundle
 gerbe $\tilde{H}_1\to H_1\toto H_0$, then
 $[\eta_G]$ maps to the Dixmier-Douady class of
the central extension $\tilde{H}_1\rtimes G\to H_1\rtimes G\toto H_0$,
under the  isomorphism $H^3_G (M)\iso  H^3_{\text{DR}} ((H\rtimes G)_\com)$
according to Theorem 4.6 in  \cite{Stienon}.
The latter corresponds exactly to $\alpha$ 
according to Corollary  \ref{pro:1.9}.

Finally, let $\theta'\in \Omega^1 (\tgm) $ be a connection
of the $S^1$-central extension $\tilde{\gm}\stackrel{\pi}{\to}\gm \toto M'$,
 which always exists according to Proposition \ref{pro:connection}.
Let $\theta =\tilde{q}^*\theta' \in \Omega^1 (\tilde{H}_1)
$, where $\tilde{q}$ is as in Eq. \eqref{eq:morphism}. Since $\tilde{q}$
is a morphism of  $S^1$-central extensions, $\theta$ is
clearly a connection of the $S^1$-central  extension
$\tilde{H}_1\stackrel{\phi}{\to} H_1\toto H_0$.
Since $\tilde{q} (\tau \cdot g)=\tilde{q} (\tau )$, $\forall \tau \in
\tilde{H}_1, \  g\in G$, it follows that $\theta$ is $G$-invariant.
Since $\hat{X}\per \theta=0$, $\forall X\in\frakg$, thus
$\theta$ is $G$-basic. This concludes the proof of the theorem.

\begin{rmk}
When  $G$ is a compact simple Lie group and
$G$ acts on $G$ by conjugation, there is an  explicit
construction of $G$-equivariant bundle gerbe due to Meinrenken
\cite{Meinrenken:gerbe} (see also \cite{GN} for the case of  $G=SU(n)$).
\end{rmk}

\subsection{Geometric transgression} 

For a Lie groupoid $\gm \toto \gm_0$,
  by $S\gm$
 we denote  the space of closed loops $\{g\in \gm |s(g)=t(g)\}$.
 Then $\gm$ acts on $S\gm$ by conjugation: 
$\gamma\cdot \tau=\gamma \tau \gamma^{-1}$,  $\forall \gamma \in \gm$ and
$\tau \in  S\gm$ such that $s(\gamma)=t(\tau )$.
 One  forms the transformation
groupoid $\Lambda\gm : S\gm\rtimes \gm \toto S\gm$, which  is called
 the {\em inertia groupoid}. 
If $\tgm\to\gm\toto \gm_0$ is an $S^1$-central extension,
then the restriction $\tgm|_{ S\gm}$
is naturally endowed with an action of
$\gm$. To see this,  for any $\gamma\in \gm$, let  $\tilde{\gamma } \in  \tgm$ 
be any of its  lifting.
Then for any  $\tilde{\tau } \in \gm|_{S\Gamma }$ such that $s(\gamma )=
t(\tilde{\tau})$, set
\begin{equation}
\label{eqn:tildeg} 
 \tilde{\tau}\cdot \gamma =\tilde{\gamma}^{-1}\tilde{\tau }{\tilde{\gamma}}.
\end{equation}
It is simple to see that this $\gm$-action is well defined, i.e.
the right hand side of Eq. \eqref{eqn:tildeg} is
independent of the choice of the lifting $\tilde{\gamma}$.
Thus $\tilde{\Gamma}|_{S\Gamma}\to S\gm$ naturally  carries 
  an $S^1$-bundle
structure over the inertia groupoid $\Lambda\gm \toto S\gm$
(see also Proposition 2.9 in \cite{TX}).

\begin{prop}\label{prop:line-bundle}
Let  $\gm\toto \gm_0$ be  a Lie groupoid.
Then any  $S^1$-central extension $\tgm\to\gm\toto \gm_0$
 induces an $S^1$-bundle
over the inertia groupoid $\Lambda\gm \toto S\gm$.
\end{prop}

\begin{rmk}
In general,  the inertia groupoid is not a Lie groupoid since
$S\gm$ may not be  a smooth manifold.
\end{rmk}

 Now assume that
$\tilde{\gm}\stackrel{\pi}{\to}\gm \toto \gm_0$ 
is  an $S^1$-central extension representing
  $\alpha \in H^3_G (M, \zz)$.
By Proposition \ref{pro:connection}, this 
central extension admits a connection $\theta \in
\Omega^1({\gmmt})$.  
According to Proposition 3.9 \cite{TX},
there exists an induced connection on the
associated $S^1$-bundle $\tilde{\gm}|_{S\gm}\to S\gm$ over the
inertia groupoid $\Lambda \gm$. 
Since $\gm \toto \gm_0$ is Morita
equivalent to $M\rtimes G\toto M$, $\Lambda \gm$
is Morita equivalent to $\Lambda (M\rtimes G)$, where
the Morita equivalence bimodule is induced by the
equivalence bimodule between $\gm\toto \gm_0$ and  $M\rtimes G\toto M$.
Thus one obtains an $S^1$-bundle $P$ with a connection
over the groupoid $\Lambda (M\rtimes G)$, according to
 Corollary 3.15 in \cite{TXL:04}.

It is clear that the groupoid 
 $\Lambda (M\rtimes G)$ is isomorphic
to $(\coprod_{g\in G}M^g ) \rtimes G\toto \coprod_{g\in G}M^g$, where
$M^g=\{x\in M|x\cdot g=x\}$ is the fixed point set under the
diffeomorphism $x\to x\cdot g$.
By  $\coprod_{g\in G} P^g \to \coprod_{g\in G}M^g$,
we denote this   $S^1$-bundle. As a consequence,
 we obtain a family
 (i.e. over a manifold instead of
over a groupoid) of $S^1$-bundles $P^g\to M^g$,
with connections, indexed by $g\in G$, on which $G$ acts
 equivariantly preserving the connections.
The main result of this section is the following

\begin{them}
\label{thm:transgression}
Let $\alpha \in H^3_G(M, \zz)$. 

\begin{enumerate}
\item 
 Assume that $\tilde{\gm}\stackrel{\pi}{\to}\gm \toto \gm_0$
is an $S^1$-central extension representing 
 $\alpha$,  
where $\gm \toto \gm_0$ is the pullback groupoid of the transformation groupoid
 $M\rtimes G\toto M$ via a surjective submersion $p: \gm_0\to M$. Then 
$\tilde{\gm}\stackrel{\pi}{\to}\gm \toto \gm_0$
 canonically induces a family of  $G$-equivariant
 flat  $S^1$-bundles $P^g\to M^g$ indexed by $g\in G$.

 Here, for any $h\in G$,  the action by $h$,
denoted  $R_h$ by abuse of notations,  is
an isomorphism of  the flat bundles $(P^g\to M^g)
 \mapsto (P^{h^{-1}gh}\to M^{h^{-1}gh})
$ over the  map $R_h :M^g \to M^{h^{-1}gh}$;
Moreover, $g$ acts on $P^g$ as an identity.

\item If  $\tilde{\gm}\stackrel{\pi}{\to}\gm \toto \gm_0$ 
and $\tilde{\gm}'\stackrel{\pi'}{\to}\gm' \toto \gm'_0$
are any two such  $S^1$-central extensions, their induced
families of  flat  $S^1$-bundles $\coprod_{g\in G}
(P^g\to M^g)$ and $\coprod_{g\in G} (P^{'g}\to M^g )$
are $G$-equivariantly isomorphic.
\end{enumerate}
\end{them}

We need a few lemmas.
Note that  $\gm\cong \{(x,g,y)\in \gm_0\times G\times \gm_0\vert\; p(x)g=p(y)\}$
and $S\gm =\{(x,g, x)\in \gm_0\times G\times \gm_0\vert\; p(x)g=p(x)\}$.
For any $g\in G$, denote 
$$(S\gm )_g=\{(x,g, x)\in \gm_0\times G\times \gm_0\vert\;
p(x)\in M^g\}.$$
By $i_g: (S\gm )_g\to \gm$, we denote the natural inclusion.

\begin{lem}
\label{lem:1.12}
Assume that $\omega\in \Omega^k (\gm )$ is a multiplicative
$k$-form, i.e. satisfies $\partial \omega=0$. Then
$i_g^* \omega =0$.
\end{lem}
\begin{pf}
For simplicity, we prove the lemma for the case $k=1$.
Differential forms of higher degree
 can be proved in a  similar manner.  Fix any tangent
vector $\delta_x \in T_x p^{-1}(M^g)$.
Let $\delta_m=p_* \delta_x\in T_m M$, where $m=p(x)$,
and $G^{\delta_m}=\{h\in G|R_{h*} \delta_m =\delta_m \}$, 
the isotropy group at $\delta_m$ of the lifted $G$-action $R_{h_*}$ 
on $TM$.  For any $h\in G^{\delta_m}$, $(\delta_x, 0_h, \delta_x)$
is clearly a well defined tangent vector in $T_{(x, h, x)} \gm$,
where $0_h\in T_h G$ is the zero tangent vector. Since
$G^{\delta_m}$ is compact, then $(\delta_x, 0_h, \delta_x)
\per \omega$, considered as a function on $G^{\delta_m}$,
must be  bounded.
On the other hand, it is simple to see that,
with respect to the tangent groupoid multiplication
$T\gm\toto T\gm_0$, we have 
$$(\delta_x, 0_g, \delta_x )\cdot (\delta_x, 0_g, \delta_x )
=(\delta_x, 0_{g^2}, \delta_x ).$$
Since $\partial \omega=0$, it follows that
$2 (\delta_x, 0_g, \delta_x )\per \omega=(\delta_x, 0_{g^2}, \delta_x )
\per \omega$.
Hence, for any $n\in \nn^*$, 
$$(\delta_x, 0_g, \delta_x ) \per \omega=\frac{1}{n}
(\delta_x, 0_{g^n}, \delta_x )\per \omega .$$
Since $g^n\in G^{\delta_m}$, it thus follows 
that $(\delta_x, 0_g, \delta_x ) \per \omega=0$.
Hence $ i_g^* \omega =0$.
\end{pf}

\begin{lem}
\label{lem:1.14}
The groupoid $\Lambda \gm|_{(S\gm )_g}\toto (S\gm )_g$ is Morita
equivalent to the transformation groupoid $M^g \rtimes G^g\toto M^g$.
\end{lem}
\begin{pf}
Consider
$$\xymatrix{
\Lambda \gm|_{(S\gm )_g}\ar[r]^\varphi\ar[d]\ar@<-.6ex>[d] &
 M^g \rtimes G^g\ar@<-.6ex>[d]\ar[d] \\
(S\gm )_g \ar[r]^{\varphi_0} &  M^g
}. $$
Here the map  $ \varphi :\Lambda \gm|_{(S\gm )_g}\to 
M^g \rtimes G^g$ is defined by $\varphi ((x, g, x), (x, h, z))=
(p(x), h), \ \forall ((x, g, x), (x, h, z)) \in S\gm \rtimes
\gm|_{(S\gm )_g}$,
  and $\varphi_0: (S\gm )_g \to M^g$ is $\varphi_0 (x, g, x)
=p(x)$, $\forall (x, g, x)\in (S\gm )_g $.
It is simple to check that the diagram above indeed
defines a Morita morphism.
\end{pf}



\noindent {\bf Proof of Theorem \ref{thm:transgression}}

(1). Note that for any $g\in G$, the $S^1$-bundle $P^g\to M^g$  
is induced from the $S^1$-bundle $\tilde{\gm}|_{(S\gm )_g}\to
 (S\gm )_g$ over $\Lambda \gm|_{(S\gm)_g}$.
According to Proposition \ref{pro:connection}, the $S^1$-central extension 
$\tilde{\gm}\to \gm$  admits a
 \connection $\theta\in \Omega^1 (\tilde{\gm})$.
Let $\omega\ \in \Omega^2 (\gm )$ be the curvature of
the $S^1$-bundle $\tgm\to \gm$, i.e. $d\theta =\pi^*\omega$.  
 Since $\partial \theta=0$,  it follows that  $\partial \omega=0$.
Hence Lemma 
  \ref{lem:1.12} implies that $\tilde{\gm}|_{(S\gm )_g}\to
 (S\gm )_g$ 
 must be a flat bundle over $\Lambda \gm|_{(S\gm )_g}\toto (S\gm )_g$. Therefore $P^g\to M^g$  is
a flat bundle over the transformation groupoid $M^g \rtimes G^g\toto M^g$,
 according to Corollary 3.15 in \cite{LTX}.

Assume that $\theta'\in \Omega^1 (\tilde{\gm})$ is another
\connection of the central extension $\tilde{\gm}\to \gm \toto \gm_0$.
 Then $\theta-\theta'
=\pi^* \xi$, where $\xi\in \Omega^1 (\gm )$ satisfies
the equation $\partial \xi =0$. Applying  Lemma \ref{lem:1.12} again,
we obtain that $\xi|_{(S\gm )_g}=0$.
 It thus follows that
$\theta'=\theta$ when being restricted to 
$\tilde{\gm}|_{(S\gm )_g}$.

By construction, for any $u\in M^g$, $P^g_u=\coprod_{\{x| p(x)=u\}}
\tilde{\gm}|_{(x, g, x)}/\sim$, where $\sim$ is
 the equivalence relation  between 
elements in $\tilde{\gm}|_{(x, g, x)}$ and those in 
$\tilde{\gm}|_{(y, g, y)}$ with
$p(x)=p(y)=u$   induced by the action  of the element $(x, 1, y)\in \Gamma$
as given by Eq. \eqref{eqn:tildeg}.
To prove that $g$ acts on $P^g$ by the identity map, it suffices
to show that $\gamma:=(x, g, x)$ acts on $P^g$ by the identity.
The latter is equivalent to the identity:

\begin{equation}
\label{eq:gxi}
\tilde{\gamma} \xi \tilde{\gamma}^{-1}=\xi, 
\end{equation}
for any $\tilde{\gamma} \in \tilde{\gm}|_{(x, g, x)}$ which lifts
$\gamma=(x, g, x)$, and any 
$\xi \in \tilde{\gm}|_{(x, g, x)}$. Since $\tilde{\gamma}$ and $\xi$ lie in
the same fiber $\tilde{\gm}|_{(x, g, x)}$, we may assume that
 $\xi=\lambda \tilde{\gamma}$, $\lambda \in \complex^*$.
 Eq. \eqref{eq:gxi}
follows immediately since $\tilde{\gm}\to\gm$ is an $S^1$-central extension.

(2). To prove (ii), assume that
$\tilde{\gm}'\stackrel{\pi'}{\to}\gm' \toto \gm'_0$
is another $S^1$-central extension representing $\alpha$,
where  $p': \gm'_0\to M$ is a  surjective submersion and $\gm'\toto \gm'_0 $ is the pullback groupoid of the transformation groupoid $M\rtimes G\toto M$
by $p'$. Let $\gm''_0=\gm_0\times_M \gm'_0$. By $\pr_1: \gm''_0\to \gm_0$
and $\pr_2: \gm''_0\to \gm'_0$, we denote the natural projections.
Let  $\gm''\toto \gm_0''$ be the pull back groupoid of
the transformation groupoid $M\rtimes G\toto M$ via $p'': \gm_0''\to M$.
By assumption, $\pr_1^* (\tgm\to \gm \toto \gm_0)$ and
$\pr_2^* (\tgm\to \gm \toto \gm_0)$
are Morita equivariant $S^1$-central extensions.
Since $p''=p'\smalcirc \pr_2=p\smalcirc \pr_1$, one   finds
  that $\gm''\toto \gm_0''$ is isomorphic
 to the   pullback groupoids $\pr_1^* ( \gm \toto \gm_0 )$
 and  $\pr_2^* (\gm'\toto \gm'_0 )$. 
 Therefore, we obtain two Morita equivariant $S^1$-central extensions:
 $\pr_1^*\tilde{\gm}\to \gm''\toto \gm_0''$
 and $\pr_2^*\tilde{\gm}'\to \gm''\toto \gm_0''$.
By Proposition 4.16 in \cite{BX}, there exists an $S^1$-bundle
$L\to \gm_0''$ such that 
$\pr_2^* \tilde{\gm}'\cong L\otimes \pr_1^* \tilde{\gm} \otimes L^{-1}$.
It thus follows that 
$\pr_2^* \tilde{\gm}'|_{(S\gm'')_g}\cong \pr_1^* \tilde{\gm}|_{(S\gm'')_g}$.
This implies that  $\tilde{\gm}'|_{(S\gm ')_g}\to (S\gm ')_g$, as an $S^1$-bundle over 
$\Lambda \gm'|_{(S\gm' )_g}\toto (S\gm' )_g$, is isomorphic to
 $\tilde{\gm}|_{(S\gm )_g}\to (S\gm )_g$,  as an $S^1$-bundle over 
$\Lambda \gm|_{(S\gm )_g}\toto (S\gm )_g$.
 
\begin{rmk}
Freed-Hopkins-Teleman also proved the existence of a family
of flat $S^1$-bundles over $M^g$ using a different method  \cite{FHT}.
\end{rmk}

\subsection{Localized twisted equivariant cohomology}

Let us first recall some basic constructions of
Block-Getzler  \cite{BG}.
Following \cite{BG}, by   a local equivariant differential form on $M$,
we mean  a smooth germ at $0\in \frakg$ of a smooth map from 
$\frakg$ to $\Omega^\bullet (M)$ equivariant  under the $G$-action.
Denote the space of all local equivariant differential forms
by  $\bOmega^\bullet_G(M)$, i.e.
$\bOmega^\bullet_G(M):= C^\infty_0(\frakg, \Omega^\bullet (M))^G$.
Here for a finite-dimensional vector space $V$,
$C^\infty_0 (V)$ denotes the algebra of germs at $0\in V$ of smooth
functions on $V$.
    It is clear that 
$\bOmega^\bullet_G(M)$  is $\zz/2$-graded, and is a module over the algebra
$C^\infty_0 (\frakg)^G$  of germs of
invariant smooth functions over $\frakg$.
The usual  equivariant differential $d_G=d+  \iota:
\Omega^\bullet_G(M)\to \Omega^{\bullet+1}_G(M)$
extends to a differential, denoted by the same symbol
$d_G$,
$$d_G=d+  \iota:
\bOmega^\bullet_G(M)\to \bOmega^{\bullet+1}_G(M)$$
satisfying $d_G^2=0$. Thus we obtain a 
$\zz/2$-graded chain complex
$(\bOmega^\bullet_G(M), d_G)$, which 
can be considered as a certain  completion of the Cartan model  of
the  equivariant cohomology $H^\bullet_G(M)$.

Now let $G$ act on the manifold underlying $G$ by conjugation:
$h\cdot g=  g^{-1}hg$.
For any $g \in  G$, by  $M^g$ we  denote the fixed point set of the
 diffeomorphism induced by $g$ on $M$. 
Let $G^g$  denote the centralizer of $g$:  $G^g =\{h \in G| gh=hg\}$, 
and $\frakg^g$ its Lie algebra.
 Consider the space of equivariant differential forms
$\bOmega\upcom (M, G)_g: =\bOmega_{G^g} (M^g )$, which consist of
germs at zero of smooth maps from $\frakg^g$ to 
$\Omega (M^g )$  equivariant under
$G^g$.
It is easy to see that if $\omega \in \bOmega\upcom (M, G)_g$,
$k\cdot \omega\in \bOmega\upcom (M, G)_{g\cdot k}$.
Moreover, the equivariant  differential coboundary
operators on $\{\bOmega_{G^g}\upcom (M^g ), \ g\in G\}$ are compatible with the $G$ action as well. That is
$$k \cdot \ed_{G^g} \omega =\ed_{G^{Ad_k g}}( k\cdot \omega), $$
where $ \ed_{G^g} :\bOmega\upcom (M, G)_g\to \bOmega\upcom (M, G)_g$
is the  equivariant differential on $\bOmega\upcom (M, G)_g$.
The family of cohomology groups 
$H^\bullet (\bOmega\upcom (M, G)_g,  \ed_{G^g})$
are called {\em localized equivariant cohomology}.

For any $\alpha\in H^3_G (M, \zz)$, let $\coprod_{g\in G}P^g$ be
a family of $G$-equivariant flat $S^1$-bundles
 as in Theorem 
\ref{thm:transgression}, and $L=\coprod_{g\in G} L^g$,
where $L^g=P^g\times_{S^1}\complex$, $\forall g\in G$,
are their associated  $G$-equivariant flat complex line bundles.

Choose a $G$-equivariant closed $3$-form $\eta_G\in Z^3_G (M)$
such that $[\eta_G]=\alpha$.
Following Block-Getzler \cite{BG}, we 
consider the  localized twisted equivariant cohomology as follows.
Denote by $\bOmega\upcom (M, G, L)_g$ the space $\bOmega_{G^g} (M^g, L^g )$ of germs at zero of $G^g$-equivariant
 smooth maps from $\frakg^g$ to $\Omega^\com (M^g , L^g)$.
By $\edd_{G^g}$ we  denote the twisted  equivariant differential operator
$\nablaa_g+ \iota-2\pi i \eta_G$
 on $\bOmega_{G^g} (M^g , L^g)$,
where $\nablaa_g: \bOmega_{G^g} (M^g , L^g)\to \bOmega_{G^g} (M^g , L^g)$
denotes the covariant differential  induced by  the flat connection
on the complex line bundle $L^g\to M^g$, and $\eta_G$ acts
on $\bOmega_{G^g} (M^g , L^g)$ by  taking the wedge product    with
$i_g^* \eta_G$. Here $i_g^*: \Omega_{G} (M)\to \Omega_{G^g} (M^g )$
is the restriction map. It is simple to see that $(\edd_{G^g})^2=0$,
and  $\{ \edd_{G^g}|g\in G\}$ are compatible with the $G$-action.
The family of cohomology groups
are denoted by $H^\bullet (\bOmega\upcom (M, G, L)_g,
\edd_{G^g})$, and are called {\em localized twisted equivariant cohomology}.

The proposition below justifies our definition.

\begin{prop}\label{prop:indep-coboundary}
Assume that $\eta_G^1$ and $\eta_G^2$
are equivariant closed  3-forms in $Z^3_G(M)$ such that
$\eta_G^2-\eta_G^1 =d_G B_G$, for some $B_G\in \Omega^2_G (M)$.
Then
\begin{eqnarray*}
\Phi_{B_G}:  \bOmega\upcom (M, G, L)_g &\to&  \bOmega\upcom (M, G, L)_g\\
\omega&\mapsto& \exp(2\pi i B_G )\omega
\end{eqnarray*}
defines  an isomorphism of the   cochain complexes
from $(\bOmega\upcom (M, G, L)_g, \nablaa_g+\iota-2\pi i \eta^1_G)$ to
$(\bOmega\upcom (M, G, L)_g, \nablaa_g+\iota-2\pi i \eta^2_G)$.

Moreover, if $B_G$ is a coboundary,
then $\Phi_{B_G}$ induces the identity map 
on the cohomology group $H^\bullet (\bOmega\upcom (M, G, L)_g,
\edd_{G^g}) $.
As a consequence,  there is a canonical map
\begin{equation}
\label{eqn:H2-Aut}
H^2_G(M)\to {\mathrm{Aut}}( H^\bullet (\bOmega\upcom (M, G, L)_g,
\edd_{G^g}) ).
\end{equation}
\end{prop}

\begin{pf}
The proof of the first part is straightforward.
We note that the exponential does make sense, since, 
if $B_G=\beta +f \in  \Omega^2 (M)^G \oplus
(\frakg^*\otimes \Omega^0 (M))^G $, 
 $\exp(2\pi i f)\in C^\infty(\Gg,M)$ is well-defined
and $\exp(2\pi i \beta)$ is
a finite sum: $\sum_{k\le {\mathrm{dim}}M/2} \frac{(2\pi i)^k}{k!}\beta^k$.

For the second part, assume  that $B_G=d_G \gamma$,
with $\gamma\in \Omega^1(M)^G$. Then $\eta_G^1=\eta_G^2$.
A simple  calculation
(see e.g. \cite[Prop. 4.8]{TX}) shows that,  for any cocycle $\omega
\in \bOmega\upcom (M, G, L)_g$,
we have 
$$\Phi_{B_G}\omega-\omega= (\nablaa+\iota-2\pi i \eta^1_G )(u\omega),$$
where $u=\gamma\frac{e^{2\pi i B_G }-1}{B_G}$. 
\end{pf}

As a consequence,  the  localized twisted  equivariant cohomology
$H^\bullet (\bOmega\upcom (M, G, L)_g, \edd_{G^g})  $
 depends on $L$
up to an isomorphism.

\section{The Hochschild-Kostant-Rosenberg theorem }

\subsection{Equivariant cyclic homology}
The Connes' Hochschild-Kostant-Rosenberg theorem states that if $M$ is a
compact manifold, then 
$$\tau:a_0\otimes\cdots\otimes a_k\mapsto
\frac{1}{k!}a_0da_1\cdots da_k$$
 induces an isomorphism from
$HP_i(C^\infty(M))$ to $H_{dR}^{i+2\zz}(M,\cc)$ for $i=0,1$. In fact,
$\tau$ is a chain map from the
periodic cyclic chain complex $(\PC_\bullet (C^\infty(M) ), b+\B)$ to
the de Rham complex $(\Omega^\bullet (M), d_{DR})$ \cite{Con85, Con:book}.
In general, we have the following

\begin{prop}
Let $A$ be  an  associative unital algebra  over $\cc$,
 $(\Omega\upcom,d)$  a differential graded algebra,
 $(\bar{\Omega}\upcom, d)$ a  chain  complex, 
and $\rho:A\to \Omega^0$ an algebra morphism. Assume that $\Tr:
\Omega\upcom\to\bar\Omega\upcom$ is a
chain map, called a    (super-)trace.
 Then the map
 $\tau: \PC_\bullet (A)\to \bar\Omega\upcom$ given by
 $$\tau:a_0\otimes \cdots\otimes a_k\mapsto
\frac{1}{k!} \Tr(\rho(a_0)d\rho(a_1)\cdots d\rho(a_k))$$
is a chain map from $(\PC_\bullet (A), b+\B)$ to 
 $(\bar{\Omega}\upcom, d)$.
\end{prop}

The goal of this section is to prove a
counterpart of the above result
for a  $G$-equivariant curved differential graded algebra.

Assume that $A$ is a topological algebra over $\cc$ 
 endowed with an  action of a compact Lie group $G$ by automorphisms.
Let $\tilde{A}=A\oplus \cc 1$ be its unitization considered
as a $G$-algebra, where $G$ acts on the unit $1$ trivially.
Set $CC^G_k(A)=C^\infty(G,\tilde{A}\otimes A^{\otimes k})^G$, where
$G$ acts on itself by conjugation, and $\otimes$ denotes an appropriate
topological tensor product chosen according to the situation.
There exist two differentials $b: CC_k^G(A)\to
CC_{k-1}^G(A)$ and $\B:CC_k^G(A)\to CC^G_{k+1}(A)$, defined,
respectively,  by

\begin{eqnarray*}
b(\varphi\otimes a_0\otimes\cdots\otimes a_k)(g)&=&
\sum_{i=0}^{k-1} (-1)^i \varphi(g) a_0\otimes\cdots\otimes
a_ia_{i+1}\otimes \cdots\otimes a_k\\
&&+(-1)^k \varphi(g) (g^{-1} a_k)a_0\otimes a_1\otimes\cdots\otimes
a_{k-1}\\
\B(\varphi\otimes a_0\otimes\cdots\otimes a_k)(g)&=&
\sum_{i=0}^{k}(-1)^{ki}\varphi(g) 1\otimes (g^{-1}a_{k-i+1})\otimes
\cdots\\
&&\quad\cdots\otimes (g^{-1}a_k)\otimes a_0\otimes \cdots\otimes a_{k-i},\\
&& \forall \varphi\in C^\infty (G), \ a_i\in A, i=1, \cdots , k.
\end{eqnarray*}

It is simple to check that $b^2=\B^2=b\B+\B b=0$. 
Let $\PC_i^G(A)=\oplus_{m\in\nn } CC_{i+2m}^G(A)$.
The $G$-equivariant periodic cyclic homology of $A$ is defined  to be the
homology group of the chain complex $(\PC_\bullet^G (A),b+\B)$:

$$HP_i^G(A)= H_i(\PC_\bullet ^G(A),b+\B).$$

The following  result  is due to  
  Brylinski \cite{Bry1, Bry2}.

\begin{prop}
Let $A$ be a topological  associative algebra, and $G$ a compact 
Lie group acting on $A$ by automorphisms.
 Then there is a natural isomomorphism 
$$HP_\bullet ^G(A) \cong  HP_\bullet 
(A\rtimes G),$$
where $A\rtimes G: =C^\infty (G,A)$ is the crossed product
algebra.
\end{prop}

We refer to   \cite{voigt} for the general theory
 of equivariant periodic cyclic homology.

\subsection{Traces on curved differential graded $G$-algebras}

We  use the standard notation for  the Lie derivative:
 $\lL_X a=X\cdot a =
\lim_{t\to 0}\frac{e^{tX}\cdot a-a}{t}$ for all $X\in \Gg$.

By an    $\nn$-graded $G$-vector space
 with a connection, we mean a $\nn$-graded
vector space  $(\Omega^n)_{n\in \nn}$
 equipped with  a degree preserving  $G$-action whose  infinitesimal
 $\frakg$-action  is denoted by
$\lL_X$, together with a  $G$-equivariant
linear map   $\nablaa: \Omega\upcom \to \Omega^{\bullet+1}$ of degree 1,
and a $G$-equivariant linear map $\iota:\Gg\to \Der^{-1}(\Omega\upcom)$ 
satisfying the following identities:
\begin{itemize}
\item[(i)] $\iota_X^2=0$ for all $X\in \Gg$;
\item[(ii)] $\nablaa \iota_X+\iota_X \nablaa = \lL_X$.
\end{itemize}
Here $\Der^{-1}(\Omega\upcom)$ denotes the space 
of degree $-1$ derivations on $\Omega\upcom$.
The operators $ \iota_X$ are called 
 contractions, the operators $\lL_X$  are called Lie derivatives,
and the operator $\nablaa$ is called connection.

A {\em morphism} between two $\nn$-graded $G$-spaces with connections is
a $G$-equivariant linear map of degree $0$ which  interchanges
the connections and contractions (and hence the Lie derivatives as well). 

Similarly, an   $\nn$-graded $G$-algebra with a connection
is a $\nn$-graded topological algebra, which is an 
$\nn$-graded $G$-vector space with a connection such that
  $\lL_X$,  $\iota_X$, and $\nablaa$  are all derivations of the graded
algebra.

Given an $\nn$-graded $G$-space with a connection 
$(\Omega^\bullet, \nabla)$, and $\Theta \in \Omega^2$,
set
\begin{equation}
\label{eq:etaG} 
\eta_G=(\nabla+  \iota)\Theta,
\end{equation}
 which is a   map from $\frakg $ to
$\Omega^\bullet$, i.e., $\eta_G (X)=\nabla \Theta+\iota_X \Theta, \ \forall X
\in\frakg$.
In the sequel, we write $\alphab =\nabla \Theta\in \Omega^3$ and
 $\eta_X=\iota_X \Theta\in \Omega^1$.

\begin{defn}
A curved differential graded  $G$-algebra is an $\nn$-graded $G$-algebra
 with a connection $(\Omega^\bullet, \nabla)$ such that 
$\nablaa^2=[\Theta,\cdot]$ for some $\Theta\in \Omega^2$ 
satisfying  the properties that $\forall X\in \frakg, \ \lL_X \Theta=0$,
 and $\eta_G(X)$ is central in $\Omega^\bullet$.
\end{defn}

For instance, for a graded $G$-manifold $M$, the de Rham complex
$(\Omega^\bullet (M), d_{\text{DR}})$ is clearly a (curved)
differential graded  $G$-algebra  with $\Theta=0$.
 
A module over a  $\nn$-graded $G$-algebra
with a connection $(\Omega^\bullet, \nabla)$
 is  a  
$\nn$-graded $G$-space
 with a  connection $(\bar\Omega^\bullet, d)$ such that 
$\forall \beta\in \Omega^\bullet$ and  $\omega\in\bar{\Omega}\upcom$,
\begin{itemize}
\item[(i)]
$d(\beta\omega)=\nabla \beta\cdot\omega+ (-1)^{|\beta|}\beta\cdot d\omega$;
\item[(ii)]
  $\bar\iota_X(\beta\omega)=(\bar\iota_X\beta)\omega + (-1)^{|\beta|} \beta
  \iota_X\omega$.
\end{itemize}

\begin{defn}
A trace map between  a curved differential graded  $G$-algebra 
  $(\Omega^n)_{n\in \nn}$
and a $\nn$-graded $G$-space  $(\bar\Omega^n)_{n\in\nn}$ with
a  connection,   is a morphism  $\Tr:\Omega^\bullet \to\bar\Omega^\bullet $
of  $\nn$-graded $G$-spaces
with connections  such that  the following identity holds for 
a  fixed    central element $g\in G$:
$$\Tr (\omega_1\omega_2)=(-1)^{|\omega_1|\,|\omega_2|}
\Tr ((g^{-1}\omega_2)\omega_1), \ \  \forall \omega_1, \ \omega_2\in
 \Omega\upcom  . $$
\end{defn}

\begin{lem}
Let $(\bar\Omega^\bullet , d)$ be  a module of a curved
 differential graded algebra  $(Z^\bullet, \nabla)$ such that $d^2=0$.
Assume that  $\nabla^2=[\Theta, \cdot]$ for $\Theta\in Z^2$.  Let 
$$\bar C_G\upcom(\bar\Omega) := C_0^\infty (\Gg,\bar\Omega\upcom)^G  $$
be the space of smooth germs at $0\in \frakg$ of   $G$-equivariant
  maps from $\frakg$ to $\bar\Omega\upcom$.
Then $(\bar C_G\upcom(\bar\Omega),  d+\bar\iota+ \eta_G)$
is a  chain complex, where the
coboundary operator $d+\bar\iota+ \eta_G$ is defined by
$$
(d+\bar\iota+\eta_G )(\omega)(X)=d(\omega (X))+\bar\iota_X(\omega(X))
+ \eta_G (X)\cdot \omega(X),
 \forall  \ \omega\in C_G\upcom(\bar\Omega),  \ X\in \frakg. $$
Here $\eta_G$ is defined by Eq. \eqref{eq:etaG}.
\end{lem}

\subsection{From equivariant periodic cyclic homology to equivariant cohomology}

We are now ready to state the main result of this section.

\begin{them}\label{thm:HKR}
Let  $(\Omega^\bullet , \nabla )$ be a curved differential graded 
$G$-algebra
with $\nablaa^2=[\Theta,\cdot]$,
and $({\bar{\Omega}}^\bullet, d)$ a $\nn$-graded
$G$-space with a connection.
Let $\Tr:\Omega^\bullet \to\bar\Omega^\bullet $ be a trace map. 
 Assume that $\bar\Omega\upcom$ is a module over
  the  curved differential graded $G$-algebra 
 $Z^\bullet$   generated by  $\eta_G(X)$, i.e.
by $\nabla \Theta$ and  $\bar\iota_X \Theta$,  $\forall X\in \frakg$,
 such that
$$\Tr (\beta\omega)=\beta \Tr \omega, \ \forall \beta \in Z^\bullet, \ 
\omega\in \Omega\upcom. $$
Then,  for any algebra homomorphism $\rho: A\to  \Omega^0$,
the map $\tau: \PC^G_\bullet (A)\to  \bar C_G\upcom(\bar\Omega)$
defined by

\begin{eqnarray*}
&&\tau(\varphi\otimes a_0\otimes\cdots\otimes a_k)(X)\\
&=&
\varphi(ge^X)\int_{\Delta_k} \Tr (\rho(a_0)\nablaa^{(t_1,X)}\rho(a_1)\cdots
\nablaa^{(t_k,X)}\rho(a_k)) e^{-\Theta} dt_1\cdots dt_k,
\end{eqnarray*}
$\forall \varphi\in C^\infty (G), \ a_i\in A, i=1, \cdots , k$,
 where $\Delta_k=\{(t_1,\ldots,t_k)\vert\; 0\le t_1\le\cdots\le t_k\le 1\}$ and
  $$\nablaa^{(t,X)}\beta = e^{-t\Theta}\nablaa(e^{-tX}\beta )e^{t\Theta}, \ \ 
\forall \beta  \in \Omega^\bullet  $$
is a chain map from  $(\PC_\bullet^G(A), b+\B)$ to 
$(C_G\upcom(\bar\Omega),d+\bar\iota+\eta_G )$.
\end{them}

We start with the following

\begin{lem}\label{lem:nablaa+i}
For all $\beta \in \Omega^0$,
$$-\frac{\del}{\del t}(e^{-t\Theta}(e^{-tX}\beta )e^{t\Theta})
=(\nablaa+\iota_X)\nablaa^{(t,X)}(\beta ).$$
\end{lem}
\begin{pf}
We have
\begin{eqnarray*}
-\frac{\del}{\del t}(e^{-t\Theta}(e^{-tX}\beta )e^{t\Theta})
&=&e^{-t\Theta}[\Theta,e^{-tX} \beta ]e^{t\Theta}+e^{-t\Theta}
  \lL_X (e^{-tX} \beta) e^{t\Theta}\\
&=& e^{-t\Theta}\nablaa^2(e^{-tX}\beta )e^{t\Theta}
  +e^{-t\Theta}\iota_X\nablaa(e^{-tX} \beta )e^{t\Theta},
\end{eqnarray*}
where we have used the identities: $\nablaa^2=[\Theta,\cdot]$ and
$\lL_X=\iota_X\nablaa+\nablaa\iota_X$.
Now, since  $\alphab=\nablaa\Theta$  is
central by assumption, we have

\begin{eqnarray*}
\nablaa\nablaa^{(t,X)}(\beta )&=&\nablaa(e^{-t\Theta}\nablaa(e^{-tX} \beta 
)e^{t\Theta})\\
&=&  -te^{-t\Theta}\alphab\nablaa(e^{-tX}\beta )e^{t\Theta}
+e^{-t\Theta}\nablaa^2(e^{-tX} \beta ) e^{t\Theta}\\
&&-e^{-t\Theta}\nablaa(e^{-tX}\beta )t\alphab e^{t\Theta}\\
&=&e^{-t\Theta}\nablaa^2(e^{-tX} \beta )e^{t\Theta}
-te^{-t\Theta}[\alpha,\nablaa(e^{-tX} \beta )]e^{t\Theta} \ \ \mbox{ (since
 $\alphab$ is central)}\\
&=& e^{-t\Theta}\nablaa^2(e^{-tX} \beta)e^{t\Theta}.
\end{eqnarray*}

Similarly, since $\iota_X\Theta:=\eta_X$  is central, we have
$\iota_X\nablaa^{(t,X)}(\beta)
  =e^{-t\Theta}\iota_X\nablaa(e^{-tX} \beta)e^{t\Theta}$.
The result thus follows.
\end{pf}

{\bf Proof of  Theorem~\ref{thm:HKR}}

Without loss of generality, we may assume that $A=\Omega^0$ 
and $\rho=\text{id}$.  We compute

\begin{eqnarray*}
\lefteqn{(d+\bar\iota)\tau(\varphi\otimes a_0\otimes\cdots\otimes
  a_k)(X)}\\
&=& \int_{\Delta_k}(d+\bar\iota_X) \Tr (\varphi(ge^X)a_0
  \nablaa^{(t_1,X)}a_1\cdots \nablaa^{(t_k,X)}a_ke^{-\Theta})\,
  dt_1\cdots dt_k\\
&=& \int_{\Delta_k} \Tr((\nablaa+\iota_X)(\varphi(ge^X)a_0
  \nablaa^{(t_1,X)}a_1\cdots \nablaa^{(t_k,X)}a_ke^{-\Theta}))\,
  dt_1\cdots dt_k\\
&=& \int_{\Delta_k} \varphi(ge^X) \Tr(\nablaa a_0
  \nablaa^{(t_1,X)}a_1\cdots \nablaa^{(t_k,X)}a_ke^{-\Theta})\,
  dt_1\cdots dt_k\\
&& +\sum_{i=1}^k (-1)^{i+1} \int_{\Delta_k} \varphi(ge^X) \Tr(a_0
  \nablaa^{(t_1,X)}a_1\cdots (\nablaa+\iota_X)\nablaa^{(t_i,X)}a_i\cdots
  \nablaa^{(t_k,X)}a_ke^{-\Theta})\,
  dt_1\cdots dt_k\\
&&+(-1)^k \int_{\Delta_k} \varphi(ge^X) \Tr(a_0
  \nablaa^{(t_1,X)}a_1\cdots \nablaa^{(t_k,X)}a_k(\nablaa+\iota_X)e^{-\Theta})\,
  dt_1\cdots dt_k \\
&=&I+II+ III.
\end{eqnarray*}

Next we examine each term $I, II$ and $III$ separately.
Since $(\nablaa+\iota_X)e^{-\Theta}=-e^{-\Theta}(\alpha+\eta_X)$, the
third term is equal to
\begin{eqnarray*}
III&=&\lefteqn{-(-1)^k\int_{\Delta_k} \Tr(\varphi(ge^X)a_0
 \nablaa^{(t_1,X)}a_1\cdots \nablaa^{(t_k,X)}a_ke^{-\Theta}(\alpha+\eta_X))\,
  dt_1\cdots dt_k}\\
&=&-(\alpha+\eta_X)\tau(\varphi\otimes a_0\otimes\cdots\otimes a_k)(X).
\end{eqnarray*}

Using Lemma~\ref{lem:nablaa+i} above, we see that the second term 
\begin{eqnarray*}
II=\lefteqn{\int_{(t_1,\ldots,\hat{t}_i,\ldots,t_k)\in\Delta_{k-1}}
\big( \sum_{i=1}^k (-1)^{i+1}
\Tr(\varphi(ge^X)a_0\nablaa^{(t_1,X)}a_1\ldots
  \nablaa^{(t_{i-1},X)}a_{i-1}
   \left(e^{-t_{i-1}\Theta}(e^{-t_{i-1}X}a_i)
    e^{t_{i-1}\Theta}\right)}\\
&& \quad\nablaa^{(t_{i+1},X)}a_{i+1}\cdots e^{-\Theta})
      +(-1)^i \Tr(\varphi(ge^X)a_0\nablaa^{(t_1,X)}a_1\cdots\\
&& \quad
  \nablaa^{(t_{i-1},X)}a_{i-1}
   \left(e^{-t_{i+1}\Theta}(e^{-t_{i+1}X}a_i)
    e^{t_{i+1}\Theta}\right) \nablaa^{(t_{i+1},X)}a_{i+1}\cdots e^{-\Theta})\,
      \big) dt_1\cdots\widehat{dt_i}\cdots dt_k
\end{eqnarray*}
with $t_0=0$ and $t_{k+1}=1$ by convention.

After re-indexing, we obtain

\begin{eqnarray*}
II=\lefteqn{\int_{(t_1,\ldots,t_{k-1})\in\Delta_{k-1}}
\sum_{i=1}^k \big(  (-1)^{i+1}
\Tr(\varphi(ge^X)a_0\nablaa^{(t_1,X)}a_1\ldots
  \nablaa^{(t_{i-1},X)}a_{i-1}e^{-t_{i-1}\Theta}}\\
&&\nablaa(e^{-t_{i-1}X})a_{i-1})
    (e^{-t_{i-1}X}a_i)e^{t_{i-1}\Theta}\cdots
  \nablaa^{(t_{k-1},X)}a_{k}e^{-\Theta})\,
      dt_1\cdots dt_{k-1}\\
&&+(-1)^i
\Tr(\varphi(ge^X)a_0\nablaa^{(t_1,X)}a_1\ldots
  \nablaa^{(t_{i-1},X)}a_{i-1}e^{-t_{i}\Theta}(e^{-t_i X}a_i)
  \nablaa(e^{-t_iX}a_{i+1})  e^{t_{i}\Theta} \cdots\\
&&\quad
  \nablaa^{(t_{k-1},X)}a_{k}\cdots e^{-\Theta}) \big )\,
      dt_1\cdots dt_{k-1}
\end{eqnarray*}

After replacing $i$ by $i+1$ in the first sum, and using the derivation
property for $\nablaa$, we obtain
\begin{eqnarray*}
II= \lefteqn{
\int_{\Delta_{k-1}}
\big( \Tr( \varphi(ge^X)a_0a_1\nablaa^{(t_1,X)}a_2\cdots
\nablaa^{(t_{k-1},X)}a_k e^{-\Theta} )}\\
&&+\sum_{i=1}^{k-1} (-1)^i
\Tr(\varphi(ge^X)a_0\nablaa^{(t_1,X)}a_1\cdots \nablaa^{(t_i,X)}(a_ia_{i+1})
\cdots \nablaa^{(t_{k-1},X)} a_k e^{-\Theta})\\
&&+ (-1)^k \Tr(\varphi(ge^X)a_0\nablaa^{(t_1,X)}a_1\cdots\nablaa^{(t_{k-1},X)}
a_{k-1}e^{-\Theta}(e^{-X}a_k)) \big )\,dt_1\cdots dt_{k-1}\\
&=& \sum_{i=0}^{k-1} (-1)^i\tau(\varphi\otimes
a_0\otimes\cdots\otimes a_ia_{i+1}\otimes\cdots\otimes a_k)\\
&&+(-1)^k\int_{\Delta_{k-1}} \Tr(\varphi(ge^X)(g^{-1}e^{-X}a_k)a_0\nablaa^{(t_1,X)}a_1
\cdots\nablaa^{(t_{k-1},X)}a_{k-1}) e^{-\Theta} \,dt_1\cdots dt_{k-1}\\
&=&(\tau\circ b)(\varphi\otimes a_0\otimes\cdots\otimes a_k)(X).
\end{eqnarray*}

It remains to show that the first term $I$ is equal to 
$(\tau\circ \B)(\varphi\otimes a_0\otimes\cdots \otimes a_k)(X)$.

\begin{eqnarray*}
\lefteqn{ (\tau\circ\B) (\varphi\otimes a_0\otimes\cdots\otimes a_k)(X)}\\
&=&\sum_{i=0}^k (-1)^{ki}
\int_{\Delta_{k+1}}\Tr(\varphi(ge^X)\nablaa^{(t_0,X)}(g^{-1}e^{-X}a_{k-i+1})\cdots\\
&&\cdots
\nablaa^{(t_{i-1},X)}(g^{-1}e^{-X}a_k)\nablaa^{(t_i,X)}a_0\cdots
\nablaa^{(t_k,X)}a_{k-i}e^{-\Theta})\,dt_0\cdots dt_k\\
&=&\sum_{i=0}^k \int_{\Delta_{k+1}}\Tr(\varphi(ge^X)\nablaa^{(t_i,X)}a_0\cdots\\
&&\cdots
\nablaa^{(t_k,X)}a_{k-i}\nablaa^{(1+t_0,X)}a_{k-i+1}\cdots
\nablaa^{(1+t_{i-1},X)}a_k e^{-\Theta})\,dt_0\cdots dt_k\\
&=& \sum_{i=0}^k \int_{\Delta_{k+1}}\Tr(\varphi(ge^X) e^{-t_i X}\cdot\left(
\nablaa a_0\cdots
\nablaa^{(t_k-t_i,X)}a_{k-i}\nablaa^{(1+t_0-t_i,X)}a_{k-i+1}\cdots\right.\\
&&\left. \cdots
\nablaa^{(1+t_{i-1}-t_i,X)}a_k\right) e^{-\Theta})\,dt_0\cdots dt_k ,
\end{eqnarray*}
where we used the $G$-invariance of $\nablaa$.

 Change  variables: $s'_0=t_i$, $s_1=t_{i+1}-t_i,
\ldots, s_{k-i}=t_k-t_i$, $s_{k-i+1}=1+t_0-t_i,\ldots,s_k=1+t_{i-1}-t_i$.
One immediately checks that $(t_0,\ldots,t_k)\in \Delta_{k+1}$
if and only if $(s_1,\ldots,s_k)\in \Delta_k$ and $s_{k-i}\le 1-s'_0
\le s_{k-i+1}$. Thus,
\begin{eqnarray*}
\lefteqn{(\tau\circ\B)(\varphi\otimes a_0\otimes\cdots\otimes a_k)(X)}\\
&=&\sum_{i=0}^k \int_{\Delta_{k}}\int_{s_{k-i}\le 1-s'_0\le s_{k-i+1}}
\Tr ( e^{-s'_0 X}\cdot\left( \varphi(ge^X) \nablaa a_0
\nablaa^{(s_1,X)}a_{1}\cdots \nablaa^{(s_k,X)}a_k\right)
e^{-\Theta})\,ds'_0ds_1\cdots ds_k.
\end{eqnarray*}

Since $\varphi\otimes a_0\otimes\cdots\otimes a_k$ is a $G$-invariant element
by assumption, it follows that  $\varphi(ge^X) \nablaa a_0
\nablaa^{(s_1,X)}a_{1}\cdots \nablaa^{(s_k,X)}a_k$
must be also $G$-invariant. 
Therefore  the equation above equals


$$\sum_{i=0}^k \int_{\Delta_{k}}(s_{k-i+1}-s_{k-i})
\Tr(\varphi(ge^X)
\nablaa a_0\nablaa^{(s_1,X)}a_{1}\cdots\nablaa^{(s_k,X)}a_k
e^{-\Theta})\,ds_1\cdots ds_k,$$
(with $s_0=0$ and $s_{k+1}=1$ by convention), which in turn
is equal to
$$\int_{\Delta_{k}}\Tr(\varphi(ge^X)
\nablaa a_0\nablaa^{(s_1,X)}a_{1}\cdots\nablaa^{(s_k,X)}a_k
e^{-\Theta})\,ds_1\cdots ds_k .$$

This concludes the proof of Theorem \ref{thm:HKR}.

\section{Main theorem}

\subsection{Pseudo-etale structure}

Let $M$ be a manifold with an action of a compact Lie group $G$.
In this section, we investigate a special kind of
$G$-equivariant bundle gerbes, whose properties
are needed for our  future study.

Recall that a  pseudo-etale structure \cite{Behrend, Tang} on a 
Lie groupoid $H_1\toto H_0$ is an integrable subbundle $F$ of $TH_1$ such that 

\begin{enumerate}
\item $F\toto TH_0 $ is a subgroupoid of the tangent Lie groupoid $TH_1
\toto TH_0 $;
\item for all $\gamma\in H_1$,  $s_*:F_\gamma\to T_{s(\gamma)}H_0 $ and $t_*:
F_\gamma\to T_{t(\gamma)} H_0 $ are isomorphisms;
\item $F|_{H_0} =TH_0$.
\end{enumerate}

\begin{defn}
\begin{enumerate}
\item
Let $f:N\to M$ be a surjective submersion.
 For any $x\in M$, denote by $N_x$ the fiber $f^{-1}(x)$.
A fiberwise measure is a family $(\mu_x)_{x\in M}$ of measures such that for all $x$, the support of $\mu_x$ is a subset of $N_x$.
\item
A fiberwise measure $\mu=(\mu_x)_{x\in M}$ is said to be smooth if for all
 $f\in C_c^\infty(N)$, the map
$$x\mapsto \int f\,d\mu_x$$
is smooth.
\end{enumerate}
\end{defn}

\begin{prop}
\label{pro:special-gerbe}
Let  $p:M'\to M$ be  an immersion such that
\begin{enumerate}
\item $\sigma:  M'\times G\to M$, $\sigma(x,g)=p(x)g$ is a
  surjective submersion, and
\item $\forall x, y\in M'$ and $g\in G$ satisfying $p(x)=p(y)g$,
there exists a diffeomorphism $\phi$
from a  neighborhood $U_x$ of $x$ in $M'$ to
  a neighborhood $U_y$ of $y$ in $M'$,
which is  compatible with the diffeomorphism
on $M$ induced by the action by $g$.  That is,
$$p(\phi (z))=p(z)g, \ \ \forall z\in U_x .$$
\end{enumerate}
Let  $H_0= M'\times G$, and   $H_1 =H_0\times_M H_0$. Then

\begin{enumerate}
\item there is a $G$-invariant   pseudo-etale structure on the 
Lie groupoid $H_1\toto H_0$;

\item  there exists a   $G$-invariant 
 Haar system $(\lambda^x)_{x\in H_0}$
 on the  Lie groupoid $H_1\toto H_0$.
\end{enumerate}
\end{prop}
\begin{pf}
(1) For any $((x, g), (y, h)) \in H_1$,
where $x, y\in M'$ and $g, h\in G$, set
$$F_{((x, g), (y, h))}=\{((u, l_g X), (v, l_h X))| \forall u\in T_x M', v\in
T_y M', X\in \frakg \ \mbox{ such that } R_{g*}p_* u=R_{h*}p_* v\} .$$
It is simple to see  that $F$ is  a $G$-invariant subbundle of
$TH_1$, and is a Lie subgroupid of  $TH_1\toto TH_0$.
To prove that $s_*:F_\gamma\to T_{s(\gamma)}H_0 $ is an
isomorphism, it suffices to show that,
for any $((x, g), (y, h)) \in H_0\times_M H_0$,
and any $u\in T_x M'$ and
$X\in \frakg$, there is  a unique $v\in
T_y M'$ such that $R_{g*}p_* u=R_{h*}p_* v$.
This holds due to  the fact
that 
$$R_{\gamma*}\big( p_* (T_xM' )\big) =p_* (T_y M'), \mbox{when } p(x)\gamma=p(y), \ x, y \in M', \gamma \in G,$$
which follows from Assumption (2). 
To see  that $F$ is integrable, note that for any
 $((x, g), (y, h))\in H_1$, the submanifold 
$\{(z, g\gamma), (\phi (z), h\gamma )|z\in U_x, \gamma \in G\}$
is a leaf of $F$ in $H_1$ through this point. 

(2) It is simple  to see that    Haar systems
 on the Lie groupoid $H_1\toto H_0$ are in one-one correspondence
with fiberwise smooth  measures of the map $\sigma :H_0\to M$, and
  $G$-invariant  Haar systems correspond to a $G$-invariant
 fiberwise smooth  measures.
    Such a measure always exists since $G$ is compact.
\end{pf}

\begin{numex}
Let $(U_i)$ be an open cover of $M$,
and  $p: M'=\coprod U_i\to M$  the covering map.
It is clear that the assumption in Proposition
\ref{pro:special-gerbe} is satisfied. 
\end{numex}


\subsection{Statement of the main theorem}

\label{sub:3.2}

Let $f:N\to M$ be a surjective submersion.
 For any  $x\in M$, denote by $N_x$ the fiber $f^{-1}(x)$.
 Assume  that $F$ is a horizontal distribution
for $f: N\to M$, i.e.  a subbundle
$F \subseteq TN$ satisfying  the condition that
 for all $y\in N$, $f_*:F_y\to T_{f(y)}X$ is an isomorphism.
For any vector field  $X$  on $M$, denote by
  $\tilde{X}$ its horizontal lifting on $N$,
and by $\Phi_t$ the flow of $\tilde{X}$.
 Note that the flow preserves fibers. More precisely,
 for any  $x\in M$ and any compact subset $K\subseteq N_x$,
 if $|t|$ is small enough,  then $\Phi_t$ is well defined on $K$ and 
maps $K$ to $N_y$ for some $y\in M$.

We say that a fiberwise smooth measure
 $(\mu_x)_{x\in M}$ of $f:N\to M$
is preserved by $F$   if for any
vector field $X$ on $M$, any $x\in M$, $f\in C_c^\infty(N_x)$, $t>0$, $y\in M$
 such that $\Phi_t$ is well defined on the support of $f$ and maps it to $N_y$, the equality
$$\int f\,d\mu_x=\int f\circ (\Phi^X_t)^{-1}\,d\mu_y$$
holds.

 Note that a  pseudo-etale structure  on a
Lie groupoid $H_1\toto H_0$  induces an action of
the Lie groupoid $H_1\toto H_0$ on the  vector
bundle $TH_0\to H_0$.

\begin{defn}
\label{def:nice}
Let $M$ be  a $G$-space.  An immersion  $p:M'\to M$ is     said to be \nice
if it satisfies the assumptions in Proposition \ref{pro:special-gerbe},
 and, in addition, there  exists a $G$-invariant
integrable  horizontal distribution $F'$ for the surjective submersion
$\sigma: H_0:= M'\times G\to M$ and a $G$-invariant
fiberwise smooth measure $(\mu_x)_{x\in M}$ for $\sigma: H_0\to M$ such that
\begin{enumerate}
\item $F'$ is preserved  under the action of $H_1\toto H_0$ induced
by the pseudo-etale structure $F\subset TH_1$;
\item $(\mu_x)_{x\in M}$ is preserved by $F'$.
\end{enumerate} 
\end{defn}

\begin{lem}
\label{lem:commutation-nabla-integration}
Assume that  $M$ is   a $G$-space, and   $p:M'\to M$ is  a \nice immersion.
Equip a $G$-invariant  Haar system on $H_1\toto H_0$ as
in Proposition \ref{pro:special-gerbe}.
For   any vector field $X$ on $M$, denote by $Y$ its horizontal lift
 to $H_0$ tangent to $F'$, and
by $\tilde{X}$ the  corresponding vector field on $H_1$
 tangent to $F$ such that $t_* (\tilde{X}_h)=Y_{t(h)}$,  $\forall h\in H_1$.
Then the flow of $\tilde{X}$ maps a $t$-fiber to another $t$-fiber, and preserves the Haar system.
\end{lem}
\begin{pf}
Denote by $\varphi_t$ and $\varphi'_t$ the flows 
of $\tilde{X}$ and $Y$ respectively.
Since $t_* (\tilde{X}_h)=Y_{t(h)}$ for all $h\in H_1$, we have
$\varphi_t(h)\in H_1^{\varphi'_t (x)}$ for all $h\in H_1^x$.

Moreover, since $s_* (\tilde{X}_h)=Y_{s(h)}$, the following diagram commutes:
$$
\xymatrix{
H_1^x \ar[r]^s\ar[d]^{\varphi_t} & H_0|_{\sigma(x)} \ar[d]^{\varphi'_t}\\
H_1^{\varphi'_t(x)} \ar[r]^s & H_0|_{\sigma(\varphi'_t(x))}.
}
$$
The conclusion follows easily. 
\end{pf}
\begin{lem}
\label{lem:nice}
 Let $M$ be  a $G$-space,  $(U_i)$  any open cover of $M$,
and  $p: M': =\coprod U_i\to M$  the covering map.
Then $p: M'\to M$ is a nice immersion.
\end{lem}
\begin{pf}
In this case,  $H_0=\coprod U_i\times G$.
  It is easy to see that for any $(x, g)\in U_i\times G$,
$$F'_{x, g}=\{(v, 0)|v\in T_x U_i\} \subset T_{(x, g)}H_0$$
is a $G$-invariant  integrable   horizontal distribution 
for the surjective submersion $\sigma: H_0 \to M$.
Define  a  $G$-invariant
fiberwise smooth measure $(\mu_x)_{x\in M}$ for $\sigma: H_0\to M$
as follows. For any $x\in M$, 
$$\sigma^{-1}(x)=\coprod_i \{(xg^{-1}, g)| \forall g\in G, xg^{-1}\in U_i\}$$
For $f\in C^\infty_c(\sigma^{-1}(x))$,
$$\int f d\mu_x=\sum_i \int_G f(xg^{-1}, g)d\lambda^G ,$$
where $\lambda^G$ is the right invariant Haar measure on $G$.
It is simple to check that, equipped with the above structures,
$\sigma$ is indeed a  \nice immersion.
\end{pf}
Now assume that    $\alpha \in H^3_G(M, \zz)  \cong 
H^3 ((M\rtimes G)\lcom, \zz )$ is  an equivariant integer class.
Let    $p: \coprod U_i\to M$ be an open cover of $M$ such that
   $\alpha$ is represented by an $S^1$-central extension
$\tilde{\gm}\stackrel{\pi}{\to}\gm \toto \coprod U_i$.
Then, according to Lemma \ref{lem:nice},
 $p:M': =\coprod U_i \to M$
 is a nice immersion.
Let $H_0=M' \times G$,   $H_1=H_0 \times_M H_0$,
 and $\tilde{H}_1=H_1\times_{\gm }\tgm$.
Then,  according to  Corollary  \ref{pro:1.9},  $\tilde{H}_1\to H_1\toto H_0$
is  a $G$-equivariant bundle gerbe over $M$ such that the Dixmier-Douady
class of $\tilde{H}_1\rtimes G
\stackrel{\phi }{\to} H_1\rtimes G \toto H_0$ is equal to
$\alpha$ under the natural
 isomorphism 
$$H^3((H\rtimes G)_\com, \zz)\cong
 H^3((M\rtimes G)_\com, \zz) \cong H^3_G(M, \zz). $$
Indeed, according to Theorem \ref{thm:main-bundlegerbe},
we can choose an equivariant connection $\theta\in\OO^1(\tilde{H}_1 )^G$
 and an equivariant curving $\bg\in\OOG^2(H_0 )$, whose
 equivariant $3$-curvature  $\etag\in\OOG^3(M)$
represents $\alpha$ in the Cartan model.
Write
\begin{equation}
\label{eq:BG}
\bg=B+\lambda, \ \ B\in \Omega^2(H_0) \text{ and } \lambda\in (\frakg^*\otimes
C^\infty (H_0))^G.
\end{equation}
Moreover, we  can assume that $\theta\in\OO^1(\tilde{H}_1 )^G$
is $G$-basic. Hence $pr_1^* \theta \in \Omega^1 ({\tilde{H}_1}\rtimes G)$
is a connection for the $S^1$-central extension $\tilde{H}_1\rtimes G
\stackrel{\phi }{\to} H_1\rtimes G \toto H_0$ according to Lemma 3.3
and Remark 3.5 \cite{Stienon}.
By Theorem \ref{thm:transgression}, there is a   family of  $G$-equivariant
 flat  $S^1$-bundles $P^g\to M^g$ indexed by $g\in G$ such that
 the action by $h$  is
an isomorphism of flat bundles $(P^g\to M^g)
 \mapsto (P^{h^{-1}gh}\to M^{h^{-1}gh})
$ over the  map $R_h :M^g \to M^{h^{-1}gh}$ and
 $g$ acts on $P^g$ as an identity map. 

Let $L^g= P^g\times_{S^1} \complex $ be the associated
complex line bundle. Denote, by $\nablaa_g$,  the induced covariant
derivative  on  $L^g\to M^g$.

By choosing a $G$-invariant Haar system $(\lambda^x)_{x\in H_0}$
 on the groupoid $H_1\toto H_0$, which always exists according to
Proposition \ref{pro:special-gerbe}, we construct a convolution
algebra $C_c^\infty(H,L)$: it  is the space of compact supported sections
of the complex line bundle $L\to H_1$, where $L=
\tilde{H}_1\times_{S^1} \complex$, equipped with the convolution
product:
$$(\xi * \eta)(\gamma)=\int_{h\in H_1^{t(\gamma)}}\xi(h)\cdot\eta
(h^{-1}\gamma)\,\lambda^{t(\gamma)}(dh), \ \ \forall \xi, \eta \in
\Gamma_c(L).$$

The following lemma  can be easily verified.

\begin{lem}
$C_c^\infty(H,L)$ is a $G$-algebra.
\end{lem}

The main result of this section is the following

\begin{them}
\label{thm:main}
Under the hypothesis above,  there exists a family of $G$-equivariant
chain maps, indexed by   $g\in G$:
$$\tau_g:(\PC_\bullet^G(C_c^\infty(H,L)),b+\B)
\to  \barom{g}
$$
By  $G$-equivariant chain maps, we mean that
the following diagram
of chain maps commutes:
$$\xymatrix{
(\PC_\bullet^G(C_c^\infty(H,L)),b+\B) \ar[r]^{\tau_{g}}
\ar[dr]_{\tau_{h^{-1}gh}} & ( \barom{g})
\ar[d]^{\phi_{h}}\\
& ( \barom{h^{-1}gh}).
}$$
\end{them}

As an immediate consequence, we have the following

\begin{them}
Under the same hypothesis as in Theorem \ref{thm:main}, there is a
  family of 
  morphisms on the level of cohomology:
\begin{equation}
\label{eq:conjecture}
\tau_g: HP^G_\bullet (C_c^\infty(H,L)) \to   H^\bullet (  \barom{g})
\end{equation}
\end{them}

The proof of Theorem \ref{thm:main} occupies the next two subsections.
The  idea is to apply Theorem \ref{thm:HKR}.
First of all, since,  $G^g$, $\forall g\in G$,
 is a Lie subgroup of $G$,  there is a natural
chain map 

\begin{equation}
(\PC_\bullet^G(C_c^\infty(H,L)),b+\B) 
\to (\PC_\bullet^{G^g} (C_c^\infty(H,L)),b+\B).
\end{equation}

To define  the map $\tau_g$,  we need to construct a chain map:
\begin{equation}
(\PC_\bullet^{G^g} (C_c^\infty(H,L)),b+\B) \to  \barom{g}
\end{equation}

 For this purpose, we need, for each fixed  $g\in G$, 

\begin{enumerate}
\item to construct a curved differential graded
$G^g$-algebra  $(\Omega^\bullet, \nabla)$ together with an algebra
homomorphism $\rho: C_c^\infty(H,L)\to \Omega^0$;
\item a trace map  $\Tr_g: \Omega^\bullet\to \Omega^\bullet (M^g,L^{g})$.
\end{enumerate}

Then we need to prove that they satisfy all the compatibility
conditions so that we can   apply Theorem \ref{thm:HKR} to  obtain  the
desired  chain map.

\subsection{A curved differential  graded $G$-algebra}

In this subsection, we deal with  the first issue
as pointed out at the end of last subsection.
After replacing $G$ by $G^g$, we may assume that $g$ is central
without loss of generality.

Denote by $\tilde\Omega^\bullet (H,L)$  the space of smooth
sections of the bundle $\wedge^\bullet F^*\otimes L \to H_1$, 
and by $\tilde\Omega_c^\bullet (H,L)$  the subspace
 of compactly supported smooth
sections.
Consider  the convolution product on $\tilde\Omega_c^\bullet (H,L)$:

$$(\xi * \eta)(\gamma)=\int_{h\in H_1^{t(\gamma)}}\xi(h)\cdot\eta
(h^{-1}\gamma)\,\lambda^{t(\gamma)}(dh).$$
Here, the product of $\xi(h)\in \wedge^k F_h^*\otimes L_h$
with  $\eta(h^{-1}\gamma)\in \wedge^l F_{h^{-1}\gamma}^*
\otimes L_{h^{-1}\gamma}$ is obtained as follows.
The product $L_h\otimes L_{h^{-1}\gamma}\to L_\gamma$ is
 induced from the groupoid structure on $\tilde{H}_1\toto H_0$, and 
$\wedge^k F_h^*\otimes \wedge^l F^*_{h^{-1}\gamma}
\to  \wedge^{k+l} F_\gamma^*$ is the composition of
the identifications $F^*_h\to F^*_\gamma$
and $F^*_{h^{-1}\gamma} \to F_\gamma^*$, via the
pseudo-etale  structure, with the wedge product.
Then $\tilde\Omega^\bullet_c(H,L)$ is a  $G$-equivariant
 graded associative  algebra. 
 The reason that we  need compactly supported ``forms'' 
on $H$ is because the convolution product does not 
generally make sense on $\tilde{\Omega}\upcom(H,L)$.

Let
$$\Omega^\com=\Omega^\com (H_0 )\oplus \tilde\Omega^\com_c(H,L).$$


 Introduce a multiplication on  $\Omega\upcom$
by the following formula:

$$(\beta_1\oplus\omega_1)\cdot (\beta_2\oplus\omega_2)
=(\beta_1\wedge \beta_2)\oplus
(t^*\beta_1\wedge   \omega_2+\omega_1\wedge   s^*\beta_2+
\omega_1 * \omega_2), \ \ \forall \beta_1, \beta_2 \in \Omega^\bullet (H_0 ),
\omega_1, \omega_2\in \tilde\Omega^\bullet_c (H,L).  $$
Note that $t^*\beta_1$ and $s^*\beta_2$ can be considered
as elements in $\tilde\Omega\upcom(H_1):=C^\infty(H_1,\wedge\upcom F^*)$,
of which  $\tilde\Omega^\bullet_c(H,L)$ is naturally a bimodule.

\begin{lem}
Under the product above, $\Omega^\bullet$ is a $G$-equivariant graded
associative algebra.
\end{lem}
\begin{pf}
It is straightforward to check the following
 relations,
$\forall \omega_i\in \tilde\Omega_c\upcom(H,L)$,  $\beta\in 
\tilde\Omega\upcom(H_1)$:

\begin{eqnarray}
t^*\beta\wedge (\omega_1*\omega_2) &=& (t^*\beta\wedge\omega_1)*\omega_2 \label{eq:1}\\
(\omega_1*\omega_2)\wedge s^*\beta &=& \omega_1*(\omega_2\wedge s^*\beta)\\
(\omega_1\wedge s^*\beta)*\omega_2 &=& \omega_1 * (t^*\beta\wedge\omega_2).
\end{eqnarray}

The associativity thus follows immediately.
\end{pf}


\begin{lem}
\label{lem:3.7}
The covariant derivative
 $\nablaa:\Omega^\com_c(H,L)\to\Omega^{\com+1}_c(H,L)$
induces an algebra derivation:
$\nablaa: \tilde\Omega^\com_c(H,L)\to\tilde\Omega^{\com+1}_c(H,L)$.
I.e.,
\begin{equation}
\label{eqn:nablaa-derivation}
\nablaa(\omega_1*\omega_2)=\nablaa\omega_1*\omega_2
+(-1)^{|\omega_1|}\omega_1*\nablaa\omega_2 \ \  \forall \omega_1, \omega_2
\in \tilde\Omega^\bullet_c(H,L).
\end{equation}
\end{lem}
\begin{pf}
Since  the pseudo-etale structure $F\subseteq TH$ is
an integrable distribution, by restriction,
 the connection on $L\to H$ induces an
$F$-connection, as a Lie algebroid connection,
 on $L$.  Therefore  we have the induced 
covariant derivative 
$\nablaa: \tilde\Omega^\com_c(H,L)\to\tilde\Omega^{\com+1}_c(H,L)$.

We now prove Eq.  (\ref{eqn:nablaa-derivation}) by induction
on $n:=|\omega_1|+|\omega_2|$. If $n=0$,  this is automatically true
since $\nablaa$ is a connection. Assume (\ref{eqn:nablaa-derivation})
is valid for $n=k$. Since $t_*: F_x \to T_{t(x)} H_0$ is an isomorphism,
 we can assume that locally
$\omega_1=t^*(df)\wedge  \eta$, where $f\in C^{\infty}(H_0 )$ and
$\eta\in \Omega^{|\omega_1|-1}_c(H,L)$. Using Eq.  \eqref{eq:1} and
 the induction hypothesis, we have

\begin{eqnarray*}
\nablaa(\omega_1*\omega_2) &=&
\nablaa(t^*(df)\wedge (\eta*\omega_2)) \\
&=& -t^*(df)\wedge\nablaa(\eta*\omega_2)\\
&=& -t^*(df)\wedge (\nablaa\eta*\omega_2+(-1)^{|\eta|} \eta *\nablaa\omega_2)\\
&=&\nablaa(t^*(df)\wedge\eta)*\omega_2+(-1)^{|\omega_1|}
((t^*df)\wedge\eta)*\nablaa\omega_2.
\end{eqnarray*}

Hence Eq. (\ref{eqn:nablaa-derivation}) follows.
\end{pf}

Extend $\nablaa$ to $\Omega\upcom$ by

\begin{equation}
\label{eq:nabla}
\nablaa(\beta\oplus\omega)=d\beta\oplus \nablaa\omega, \ \
\forall \omega\in \tilde\Omega_c\upcom(H,L), \ \ \beta\in \Omega\upcom(H_0).
\end{equation}

The following lemma is straightforward:

\begin{lem}
$\nablaa: \Omega\upcom  \to \Omega^{\bullet+1}$ is a degree $+1$
 derivation.
\end{lem}

\begin{lem} 
The map $\nablaa: \Omega\upcom  \to \Omega^{\bullet+1}$ satisfies
$$\nablaa^2=-2\pi i [B,\cdot],  $$ 
where $B\in \Omega^2(H_0 )\subset \Omega\upcom$ is defined as
in Eq. \eqref{eq:BG}.
\end{lem}

\begin{pf}
The relation $\nablaa^2\beta=-2\pi i [B,\beta ]$ holds automatically
 for $\beta \in \Omega\upcom(H_0 )$ since both sides are zero. 
For  $\omega\in \tilde\Omega_c\upcom(H,L)$, we have 

\begin{eqnarray*}
\nablaa^2\omega&=&
2\pi i d\theta\wedge  \omega = 2\pi i \del B\wedge \omega
=-2\pi i (t^*B\wedge \omega-s^*B\wedge  \omega)= -2\pi i [B,\omega].\\
\end{eqnarray*}
The conclusion thus follows.
\end{pf}

\begin{prop}
$(\Omega\upcom, \nablaa)$ is a curved differential graded  $G$-algebra.
\end{prop}
\begin{pf}
Define the map $\iota_X$, $\forall X\in \frakg$,
 on $\Omega\upcom(H_0 )$ by the usual contraction.
Since the  pseudo-etale structure $F\subset TH_1$ is preserved by the 
group $G$-action, it follows that the contraction
 map $\iota_X$ is also defined on
${\Omega}\upcom_c(H,L)$. Thus we have a map
$\iota_X:\Omega\upcom\to \Omega^{\bullet-1}$ satisfying
$\iota_X^2=0$. One proves,  by induction similar to 
 the proof of Lemma \ref{lem:3.7}, that $\iota_X$ is indeed a 
derivation.

It remains to prove that  $\nablaa\iota_X+\iota_X\nablaa=\lL_X$ 
holds. This  is clearly true on the component 
$\Omega\upcom(H_0)$.  Since both sides are derivations
on $\tilde\Omega^\bullet _c(H,L)$ and this holds for
elements in  $\tilde\Omega\upcom(H):=C^\infty(H_1 ,\wedge\upcom F^*)$,
it suffices to prove this identity for elements of degree
zero, i.e.  $\lL_X \xi=\nablaa_X\xi$ for all $X\in \Gg$ and 
$\xi\in C_c^\infty(H_1 , L)$.
Note that  for a $G$-equivariant complex line bundle we always have
the identity $\nablaa_X -\lL_X=\iota_X\theta$.  Here $\iota_X\theta$
can be considered as a function on $H_0$.

Since $\theta$ is $G$-basic, it follows that $\lL_X \xi=\nablaa_X\xi$.
This concludes the proof.
\end{pf}

\subsection{Trace map}
Now  we deal with  the second  issue
 at the end of Section \ref{sub:3.2}.

Consider  $\bar\Omega^n=\Omega_c^n(M^g,L^{g})$.
Let $d: \bar\Omega^n\to \bar\Omega^{n+1}$ denote the  covariant
derivative $\Omega_c^n(M^g,L^{g})\to \Omega_c^{n+1}(M^g,L^{g})$ 
of the   flat complex line bundle $L^{g}\to M^g$.
It is clear that $\bar\Omega^\bullet$ admits
a $G^g$-action. Let $\bar{\iota}: \bar{\Omega}^n\to \bar{\Omega}^{n-1}$
be the usual contraction.

\begin{lem}
$(\bar\Omega , d)$ is a     $\nn$-graded $G$-vector space
 with a connection.
\end{lem}

For any fixed $g\in G$, since $G^g$ is a Lie subgroup of $G$,
$(\Omega^\bullet, \nabla)$ is a  curved $G^g$-differential
algebra.  Next we will  introduce  a trace map
 $$\Tr_g:\Omega^\bullet \to \bar\Omega^\bullet $$

On the first component
$\Omega^\bullet (H_0 )$, we  set   $\Tr_g$ to be zero.
Before we define $\Tr_g$ on the second component of $\Omega\upcom$,
some remarks are in order.
Denote $(S(H_1\rtimes G))_g$ by $H_1^g$. I.e.
$$H_1^g =\{(x, xg^{-1})|\sigma (x)\in M^g\}. $$

Note that the subbundle $F'\subset TH_0$ corresponds to
a subbundle $F_1\subset F$ such that $t_*(F_1)_\gamma
=F'_{t(\gamma)}$ for all $\gamma\in H_1$. Let $F^{g}_1=
F_1\cap TH_1^{g}$. Then $(\sigma\smalcirc s)_*=(\sigma\smalcirc t)_*:
(F_1)_\gamma\to TM^g_{(\sigma\smalcirc s)(\gamma)}$ is an isomorphism for all
$\gamma\in H_1^{g}$.
Denote by $i^{g}: \wedge F^*\to \wedge (F_1^{g})^*$ the restriction map.

We define, for any fixed  $\gamma$,

$$(\Tr_g\omega)_\gamma=i^{g}
\int_{h\in H_1^x}h_*\omega_{h^{-1}\gamma h^{g^{-1}}}\,\lambda^x(dh)
\in\wedge (F_1^{g})_\gamma^*\otimes L|_\gamma .$$

If we write $\gamma=(x,xg^{-1}) , \ x\in H_0^{g}$, and $h=(x, y)$,
then 
$$h^{-1}\gamma h^{g^{-1}}=(y, x)(x, xg^{-1}) (xg^{-1}, yg^{-1})
=(y, yg^{-1}). $$
 Therefore $ \omega_{h^{-1}\gamma h^{g^{-1}}}\in
\wedge F^*_{(y, yg^{-1})}\otimes L|_{(y, yg^{-1})}$  and
$h_*\omega_{h^{-1}\gamma h^{g^{-1}}}\in \wedge
 F^*_{(x, xg^{-1})}\otimes L|_{(x, xg^{-1})}$. Hence under the map
 $i^{g}$, it goes to
$\wedge (F_1^{g})_{(x, xg^{-1})}^*\otimes L|_{(x, xg^{-1})}$. 
Since  $\Tr_g \omega$ is $H_1^{g}$-invariant,  it
 defines an element in  $\Omega\upcom(M^g,L^{g})$.

\begin{lem}
We have
$$\Tr_g(\omega_1*\omega_2)= (-1)^{|\omega_1||\omega_2|} \Tr_g((g^{-1}\omega_2)*\omega_1)$$
\end{lem}
\begin{pf}
We compute
\begin{eqnarray*}
\Tr_g(\omega_1*\omega_2)_\gamma&=& i^{g} 
\int_{h\in H_1^x}\int_{h'\in H_1^{s(h)}} h_*
((\omega_1)_{h'}\cdot (\omega_2)_{(h')^{-1}h^{-1}\gamma h^{g^{-1}}})\,
\lambda^{s(h)}(dh')\lambda^x(dh)\\
&=&  i^{g}\int_{h,k\in H_1^x}h_*(
(\omega_1)_{h^{-1}k} \cdot(\omega_2)_{k^{-1}\gamma h^{g^{-1}}})\,
\lambda^{x}(dh)\lambda^x(dk)\quad(k=hh')\\
\Tr_g((g^{-1}\omega_2)*\omega_1)_\gamma&=&
 i^{g}\int_{h,k\in H_1^x}k_*(
(g^{-1}\omega_2)_{k^{-1}h} \cdot(\omega_1)_{h^{-1}\gamma k^{g^{-1}}})\,
\lambda^{x}(dh)\lambda^x(dk).
\end{eqnarray*}

Replacing $k$ by $\gamma k^{g^{-1}}$ in the expression for
$\Tr_g(\omega_1*\omega_2)_\gamma$, we get
$$\Tr_g(\omega_1*\omega_2)_\gamma= i^{g}\int_{h,k\in H_1^x}h_*(
(\omega_1)_{h^{-1}\gamma k^{g^{-1}}} 
\cdot(( g^{-1} \omega_2)_{k^{-1} h})^{g^{-1}})\,
\lambda^{x}(dh)\lambda^x(dk).$$

It  thus remains to prove that
$$i^{g}(h^{-1}k)_*(\alpha_2\cdot\alpha_1)= (-1)^{|\alpha_1||\alpha_2|}
i^{g}\alpha_1\cdot \alpha_2^{g^{-1}}$$
for all $\alpha_1\in(\wedge F^*\otimes L)|_{h^{-1}\gamma k^{g^{-1}}}$
and $\alpha_2\in(\wedge F^*\otimes L)|_{k^{-1}h}$.
Replacing  $\gamma$ by $k^{-1}\gamma k^{g^{-1}}$ and $h$ by $k^{-1}h$
to simplify notations, it suffices to show:
$$i^{g} [ (h^{-1})_*(\alpha_2\cdot\alpha_1)]
    = (-1)^{|\alpha_1||\alpha_2|}
i^{g}[ \alpha_1\cdot (\alpha_2^{g^{-1}} ) ]$$
for all $\alpha_1\in(\wedge F^*\otimes L)|_{h^{-1}\gamma}$
and $\alpha_2\in(\wedge F^*\otimes L)|_{h}$.
We may assume that $\alpha_1=\eta_1\otimes\xi_1$ and
$\alpha_2=\eta_2\otimes\xi_2$, for $\eta_1, \eta_2\in
\wedge F^*$ and $\xi_1, \xi_2\in L$.
 It then suffices to establish the following
equalities:
\begin{itemize}
\item[(a)] $(h^{-1})_*(\xi_2\cdot\xi_1)=\xi_1\cdot\xi_2^{g^{-1}}$;
\item[(b)] $i^{g}(h^{-1})_*  [ (r_{h^{-1}\gamma}\eta_2)\wedge
(l_h\eta_1)  ] = i^{g} (-1)^{|\eta_1|\,|\eta_2|}(r_{h^{g^{-1}}}\eta_1)
\wedge  (l_{h^{-1}\gamma}\eta_2^{g^{-1}} )$.
\end{itemize}

For (a), choose a lift $\tilde{h}\in\tilde{H_1}$ of $h$, and identify
$\tilde{h}$ to $(\tilde{h},1)\in \tilde{H_1}\times_{S^1}\cc=L$.
Then, from (\ref{eqn:action}),
$(h^{-1})_*(\xi_2\cdot\xi_1)=(\tilde{h})^{-1}(\xi_2\cdot\xi_1)
({\tilde{h}}^{g^{-1}})=({\tilde{h}}^{-1}\xi_2)(\xi_1{\tilde{h}}^{g^{-1}})$.
Now, ${\tilde{h}}^{-1}\xi_2$ is an element of $L_{s(h)}\cong \cc$,
hence can be identified to the complex number
$({\tilde{h}}^{-1}\xi_2)^{g^{-1}}\in L_{s(h)g^{-1}}\cong \cc$.
Therefore we have $(h^{-1})_*(\xi_2\cdot\xi_1)=(\xi_1{\tilde{h}}^{g^{-1}})
({\tilde{h}}^{-1}\xi_2)^{g^{-1}}=\xi_1\xi_2^{g^{-1}}$ as claimed.
\par\medskip

For (b), using the fact that $(h^{-1})_*=l_{h^{-1}}r_{h^{g^{-1}}}$,
we see that (b) reduces to
$$i^{g} l_{h^{-1}} r_{h^{-1}\gamma h^{g^{-1}}}\eta_2
=i^{g} l_{h^{-1}\gamma}\eta_2^{g^{-1}}.$$
We can of course assume that $\eta_2\in F^*_h$. By duality,
this is equivalent to
$$ l_{h^{-1}} r_{h^{-1}\gamma h^{g^{-1}}}X
= l_{h^{-1}\gamma}X^{g^{-1}},$$
i.e. to $r_{h^{-1}\gamma h^{g^{-1}}}X
= l_{\gamma}X^{g^{-1}}$
for all $X\in (F_1)_h^{g}$. Since $(\sigma\circ s)_*:
(F_1)_{\gamma h^{g^{-1}}}\to T_{\sigma\circ s(h)}M^g$ is an isomorphism,
it suffices to check that both sides coincide
after  applying the map $(\sigma\circ s)_*$.

$$(\sigma\circ s)_*(l_\gamma X^{g^{-1}})=(\sigma\circ s)_*(X^{g^{-1}})
=((\sigma\circ s)_*(X)){g^{-1}}=(\sigma\circ s)_*(X)$$
since $(\sigma\circ s)_*(X)\in M^g$ by assumption. On the other hand,

$$(\sigma\circ s)_*(r_{h^{-1}\gamma h^{g^{-1}}} X)=
(\sigma\circ t)_*(r_{h^{-1}\gamma h^{g^{-1}}} X)=
(\sigma\circ t)_*(X)= (\sigma\circ s)_*(X).$$

This completes the proof.
\end{pf}

\begin{lem}
We have
$$\Tr_g\circ\nablaa=\nablaa_g\circ \Tr_g .$$
\end{lem}
\begin{pf}
Since the map $\Tr_g$ factors through the restriction  to $M^g$,
we may simply assume that $M=M^g$, and thus $H_0=H_0^{g}$.

For all $y\in H_0$, the element $\omega_{(y,yg^{-1})}\in \Lambda F^*_{(y,yg^{-1})}\otimes L_{(y,yg^{-1})}$ restricts to an element
 $\omega'_y\in \Lambda {F'}^*_y\otimes L_{(y,yg^{-1})}$.

For all $x\in H_0$, we have
$$(\Tr_g\omega)(X_1,\ldots,X_n)(\sigma(x))
=\int_{h=(x,y)\in H_1^x} h_*\omega'(\tilde{X}_1(y),\ldots,\tilde{X}_n(y))\,\lambda^x(dh).$$

Then

\begin{eqnarray*}
 \lefteqn{(\nablaa_g(\Tr_g\omega))(X_1,\ldots,X_n)(x)}\\
&=& \sum_{i=1}^n (-1)^{i+1}\nablaa^g_{X_i}\cdot (\Tr_g\omega)(X_1,\ldots,\widehat{X_i},
\ldots,X_n)\\
&&\quad + \sum_{1\leq i<j\leq n} (-1)^{i+j} (\Tr_g\omega)([X_i,X_j],X_1,\ldots, \widehat{X_i},
\ldots,\widehat{X_j},\ldots,{X_n})\\
&=& \sum_{i=1}^n (-1)^{i+1}\nablaa^g_{X_i}\cdot h_* \int_{h\in H_1^x}\omega'(\tilde{X}_1(y),\ldots,\widehat{\tilde{X}_i(y)},
\ldots,\tilde{X}_n(y))\,\lambda^x(dh)\\
&&\quad + \sum_{1\leq i<j\leq n} (-1)^{i+j} \int_{h\in H_1^x} h_* \omega'(\widetilde{[X_i,X_j]}(y),\tilde{X}_1(y),\ldots, \widehat{\tilde{X}_i(y)},
\ldots,\widehat{\tilde{X}_j(y)},\ldots,{\tilde{X}_n(y)})\,\lambda^x(dh)\\
&=& \sum_{i=1}^n (-1)^{i+1}\nablaa^g_{X_i}\cdot \int_{h\in H_1^x} h_*\omega'(\tilde{X}_1(y),\ldots,\widehat{\tilde{X}_i(y)},
\ldots,\tilde{X}_n(y))\,\lambda^x(dh)\\
&&\quad + \sum_{1\leq i<j\leq n} (-1)^{i+j} \int_{h\in H_1^x} h_*\omega'({[\tilde{X}_i,\tilde{X}_j]}(y),\tilde{X}_1(y),\ldots, \widehat{\tilde{X}_i(y)},
\ldots,\widehat{\tilde{X}_j(y)},\ldots,{\tilde{X}_n(y)})\,\lambda^x(dh).
\end{eqnarray*}

The last equality is a consequence of integrability of the bundle $F'$.

To conclude, we need to show that in the expression $\nablaa^g_{X_i}\cdot \int_{h\in H_1^x}\omega'(\tilde{X}_1(y),\ldots,\widehat{\tilde{X}_i(y)},
\ldots,\tilde{X}_n(y))\,\lambda^x(dh)$, derivation commutes with integration, i.e. that for every $L$-valued section the equality
$$\nablaa_{\tilde{X}} \int_{h=(x,y)\in H_1^x} h_*\xi(y)\,\lambda^x(dh) = \int_{h=(x,y)\in H_1^x} h_* ((\nablaa_{\tilde{X}}\xi)(y))\,\lambda^x(dh)$$
holds; this is a consequence
of Lemma~\ref{lem:commutation-nabla-integration}.
\end{pf}

The following results can be easily verified directly.

\begin{prop}
 The family of trace maps $\Tr_g: \Omega^\bullet \to \bar\Omega_g$
is $G$-equivariant, i.e., the following diagram

$$\xymatrix{
\Omega^\bullet \ar[r]^{\Tr_g}
\ar[dr]_{\Tr_{h^{-1}gh} } & \bar\Omega_g:= \Omega_c^\bullet (M^g, L^g)
\ar[d]^{(\phi_h^{-1})^* }\\
&  \bar\Omega_{h^{-1}gh }:= \Omega_c^\bullet (M^h, L^h)
}$$
commutes.
\end{prop}

\section{Discussions and open questions}
\label{sec:discussion}

\subsection{Global equivariant differential forms a la Block-Getzler}

In this section, we briefly recall the basic construction of
global equivariant differential forms a la Block-Getzler. We closely
follow the approach of \cite{BG}.

Recall that $G$ acts on the manifold underlying $G$ by conjugation:
 $h\cdot g=  g^{-1}hg$.
Equip $G$ with the topology of invariant open sets: 
\begin{equation}
\label{eq:open}
\calo=\{U \subset  G \mbox{ open} | U =  U \cdot g \ \ \forall g\in  G \}.
\end{equation}
 Construct a sheaf $\bOmega\upcom(M, G)$ over $G$ as follows.
    The stalk of the sheaf $\bOmega\upcom (M, G)$
 at $g\in   G$ is the space of equivariant differential forms
$\bOmega\upcom (M, G)_g=\bOmega_{G^g} (M^g )$, which consist of
germs at zero of smooth maps from $\frakg^g$ to 
$\Omega (M^g )$  equivariant under
$G^g$.
It is easy to see that if $\omega \in \bOmega\upcom (M, G)_g$,
$k\cdot \omega\in \bOmega\upcom (M, G)_{g\cdot k}$.
Therefore the group $G$ acts on the sheaf
$\bOmega\upcom (M, G)$ in a way
compatible with its conjugation action on $G$.
Moreover, the equivariant  differential coboundary
operators on $\{\bOmega_{G^g}\upcom (M^g ), \ g\in G\}$ are compatible with the $G$ action as well. That is
$$k \cdot \ed_{G^g} \omega =\ed_{G^{Ad_k g}}( k\cdot \omega), $$
where $ \ed_{G^g} :\bOmega\upcom (M, G)_g\to \bOmega\upcom (M, G)_g$
is the  equivariant differential on $\bOmega\upcom (M, G)_g$.

\begin{defn}
We say that a point $h=g\exp X\in G^g$, $X\in \frakg^g$,
 is near the  point $g\in  G^g$
if $M^{h}\subseteq M^g$ and $G^h \subseteq G^g$.  
\end{defn}

Note that, from a theorem of Mostow-Palais \cite{Mostov, Palais},
it follows that the set of all points in $ G^g$ near
$g$ is indeed an open neighborhood of $g$ in $G$ under the topology
given as in
 \eqref{eq:open}.  Hence a section $\omega \in \gm (U,   \bOmega\upcom (M, G))$
 of the sheaf $\bOmega\upcom (M, G)$ over an invariant open set
$U \subset  G$ is given by, for each point $g\in
  U$, an element $\omega_g \in \bOmega\upcom (M, G)_g$, such
that if $h=g\exp{X}\in G^g$ is near $g$, we have the equality of germs:
$$\omega_h (Y)=\omega_g (X+Y) \in \bOmega\upcom (M, G)_h, \  \ \ \forall Y\in \frakg^h . $$
Thus $\bOmega\upcom (M, G)$ is an equivariant sheaf of differential
 graded algebras over $G$.
By a global equivariant differential form on $G$, we mean
an equivariant section $\omega \in \gm (G,   \bOmega\upcom (M, G))^G$,
i.e. $\omega_{k\cdot g}= k\cdot\omega_g$, 
$\forall g, k\in G$. 

\begin{rmk}
To understand the meaning of
the above conditions on global equivariant differential forms,
it is  useful to consider the following simple example in an anagous
 situation.  Let $f\in C^\infty(G)^G$. For each fixed $g\in G$,
denote $f_g(X)$ the germ of the function $X\to f(g\exp{X})$, $\forall
X\in \frakg^g$. It is easy to check that the following identities hold
\begin{enumerate}
\item if $h=g\exp{X}\in G^g$, $X\in \frakg^g$, then $f_h (Y)=f_g (X+Y)$,
 $\forall Y\in \frakg^g$;
\item for any $r\in G$,  $f_{r^{-1}gr} =({Ad_{r^{-1}}})^* f_g$.
\end{enumerate}
\end{rmk}

Let 
$$\cala_G\upcom(M) =\gm (G,  \bOmega\upcom (M, G))^G$$
be the space of  global equivariant differential forms on $G$.
The family of differentials $\{\ed_{G^g}\}_{g\in G}$ induces a differential
$d_{\text{eq}}$ on $\cala_G\upcom(M)$. The 
delocalized  equivariant cohomology
$H_{G, delocalized}\upcom (M)$ is defined as its
$\zz/2$-graded cohomology: 
$$H_{G, delocalized}\upcom (M):=H\upcom (\cala_G\upcom(M), d_{\text{eq}}).$$

The following result is due to Block-Getzler \cite{BG}   when
the  Lie group $G$ is compact,
 and to Baum-Brylinski-MacPherson \cite{BBM} when $G=S^1$. 

\begin{them}
Let $G$ be a compact Lie group, $M$ a compact manifold on which $G$
 acts smoothly. Then
$$HP_\com^G (C^\infty (M))\cong H\upcom (\cala_G\upcom(M), d_{\text{eq}}).$$
\end{them}

\subsection{Delocalized twisted equivariant cohomology}

A natural open question arises:

\begin{question}
Introduce global twisted equivariant
differential forms by modifying   the notion of global equivariant
differential forms of Block-Getzler \cite{BG} in the previous section,
 and  define  delocalized twisted equivariant cohomology.
\end{question}

For any $\alpha\in H^3_G (M, \zz)$, let $\coprod_{g\in G}P^g$ be
a family of $G$-equivariant flat $S^1$-bundles
 as in Theorem 
\ref{thm:transgression}, and $L=\coprod_{g\in G} L^g$,
where $L^g=P^g\times_{S^1}\complex$, $\forall g\in G$,
are their associated  $G$-equivariant flat complex line bundles.

Following Block-Getzler \cite{BG}, the key issue is to introduce 
a sheaf $\bOmega\upcom (M, G, L)$ over $G$  whose
  stalk
 at $g\in   G$ is the space of equivariant differential forms, with
coefficients in $L^g$:
$\bOmega\upcom (M, G, L)_g=\bOmega_{G^g} (M^g, L^g )$. 

For this purpose, we  propose the following

\begin{conjecture}
Assume that  a point $h=g\exp X\in G^g$, $X\in \frakg^g$,
 is near the  point $g\in  G^g$. Then there exists
a canonical  isomorphism of flat $S^1$-bundles over $M^h$:
\begin{equation}
\label{eq:phigh}
\phi_{gh}: P^h\to P^g|_{M^{h}},
\end{equation}
where $P^g|_{M^{h}}\to M^h$ denotes
 the restriction of $P^g\to M^g$ under the inclusion $M^{h}\subseteq M^g$.  
\end{conjecture}

If the conjecture above holds, one can make sense of 
{\em a global twisted  equivariant differential form}
 on $G$ by defining it to be
an  equivariant section $\omega \in \gm (G,   \bOmega\upcom (M, G, L))^G$,
i.e. a family $\{\omega_g\in \bOmega_{G^g} (M^g , L^g)|g\in G\}$ such that

\begin{enumerate}
\item  $\omega_{k\cdot g}= k\cdot\omega_g$,
$\forall g, k\in G$,   and
\item  if $h=g\exp X\in G^g, \ X\in \frakg^g$ is near $g$, we have the equality of germs
$$ \phi_{gh} [\omega_h (Y)]=\omega_g (X+Y)|_{M^h}, 
\in \bOmega\upcom (M, G, L)_h, \  \ \
\forall Y\in \frakg^h , $$
where $\phi_{gh} : L^h\to L^g$ is the  canonical isomorphism
induced by the isomorphism  of $S^1$-bundles as in Eq. \eqref{eq:phigh}.
\end{enumerate}

Let $\cala_G\upcom(M, L) =\gm (G,  \bOmega\upcom (M, G, L))$
be the space of global twisted  equivariant differential forms on $G$.
The family of differentials $\{\edd_{G^g}\}_{g\in G}$
 induces a differential
$d_{\text{eq}}^{\alpha}$ on $\cala_G\upcom(M, L)$.
By abuse of notations, we write
$$d_{\text{eq}}^{\alpha}=
\nablaa+\iota-2\pi i \eta_G. $$
	 The delocalized twisted equivariant cohomology 
$H_{G, delocalized, \alpha}\upcom (M)$ can then  be defined as its
$\zz/2$-graded cohomology:
$$H_{G, delocalized, \alpha}\upcom (M):=H\upcom (\cala_G\upcom(M, L),
 d_{\text{eq}}^{\alpha}).$$

\begin{rmk}
A  possible way to construct $\phi_{gh}$
in Eq. \eqref{eq:phigh} is to use the parallel transportation along the 
path $t\mapsto g\exp(tX)$. More precisely,
 assume that  $S^1\to \tilde{\Gamma}\to \Gamma \toto M'$
 is an $S^1$-central extension representing
$\alpha \in H^3_G(M, \zz)$,
 where $f:M'\to M$ is a surjective submersion.  Let $\nabla
\in \Omega^1 (\tilde{\Gamma})$  be a gerbe connection.
The fibers $P^g_x$ and $P^h_x$ can be identified with $P_{x,g,x}$ and 
$P_{x,h,x}$, respectively.
 Then parallel transportation along the path $t\mapsto (x,g\exp(tX),x)$
 can thus be used to identify $P_{x,g,x}:=\tilde{\Gamma}|_{(x,g,x)}$ with 
$P_{x,h,x}:=\tilde{\Gamma}|_{(x,h,x)}$.
However, for this identification to intertwine the connections
 $\nabla^g$ and $\nabla^h$, one  needs   that
the curvature 2-form $\omega \in \Omega^2(\Gamma)$ of the
connection $\nabla$ satisfies
the condition  
\begin{equation}
\label{eq:flatgh}
\omega_{x,g,x}((v,0,v),(0,X,0))=0,
\end{equation}
$\forall x\in M^h\subset M^g$ and all $v\in T_x M^h$.
This condition does not necessarily always hold as indicated
by the example below.

Let  $M$ be $S^1$ endowed with a trivial action of $G=S^1$. 
Let $\gm\toto M$ be the transformation groupoid
 $M\rtimes G \toto M$.
Consider
$$\tilde{\Gamma}=\frac{[0,1]\times G\times S^1}{(0,g,\lambda)\sim (1,g,g+\lambda)}.$$
(The product in $S^1$ is written additively.)
It is clear that $M\times S^1\to \tilde{\Gamma}\to\Gamma \toto M$
is an $S^1$-central extension, where
the map $\tilde{\Gamma}\to\Gamma$ is defined as $(u,g,\lambda)\mapsto (u,g)$
 (with $S^1$ identified to $\rr/\zz$), and the map $M\times S^1\to\tilde{\Gamma}$ is $(u,\lambda)\mapsto (u,0,\lambda)$.
Since $\tilde{\Gamma}\to\Gamma$ is a non-trivial $S^1$-principal bundle,
its first Chern class must be nonzero.  Therefore, Condition \eqref{eq:flatgh}
fails in this case.
\end{rmk}

\subsection{De Rham model of equivariant twisted K-theory}

Let $A$ be a topological  associative algebra, and $G$ a compact
Lie group acting on $A$ by automorphisms.
 Then there is an equivariant Chern character \cite{Block, Bry1, Bry2}:

\begin{equation}
\label{eq:eqchern}
\text{ch}^G: K_\bullet^G(A) \to  HP_\bullet^G(A) 
\end{equation}
from the equivariant  $K$-theory of $A$ to the periodic cyclic
homology $HP_\bullet^G(A)$.

Let $R(G)$ be  the representation ring
of $G$,  and $R^\infty (G)$ the algebra $C^\infty(G)^G$ of
smooth functions on the group $G$ invariant under the conjugation.
Since the character map sends $R(G)$   to $R^\infty (G)$,
$R^\infty (G)$ is an algebra over the ring $R(G)$.
The following result is due to Block \cite{Block}
and Brylinski \cite{Bry1, Bry2}.

\begin{them}
\label{thm:Block-Bry}
Let $G$ be a compact Lie group and  $A$  a topological  $G$-algebra.
 Then the equivariant Chern character \eqref{eq:eqchern} induces an isomorphism
\begin{equation}
HP_\bullet ^G(A) \cong K_\bullet^G(A)\otimes_{R(G)}R^\infty (G)
\end{equation}
\end{them}

Apply the theorem above to the   topological algebra: 
$$A= C_c^\infty(H,L).$$
 By definition (see \cite{TXL:04} for details), $K_\bullet^G(C_c^\infty(H,L))$
is exactly the twisted $K$-theory group $K\upcom_{G, \alpha}(M)$.
Thus we obtain

\begin{cor}
Under the same hypothesis as in Theorem \ref{thm:main}, we have
\begin{equation}
HP_\bullet ^G (C_c^\infty(H,L))\cong K\upcom_{G, \alpha}(M)\otimes_{R(G)}R^\infty (G)
\end{equation}
\end{cor}

Therefore one may think of $HP_\bullet ^G (C_c^\infty(H,L))$
as a de Rham model of equivariant twisted K-theory.


\begin{conjecture}
The family of chain maps $\{\tau_g\}_{g\in G}$ as in Theorem
\ref{thm:main} induces a quasi-isomorphism
\begin{equation}
\label{eq:tau}
  (\PC_\bullet^G(C_c^\infty(H,L)),b+\B)\to
(\cala_G\upcom(M, L), d_{\text{eq}}^{\alpha})
\end{equation}
\end{conjecture}

An immediate consequence yields the following

\begin{conjecture}
\label{conjecture2}
Let $G$ be a compact Lie group, $M$ a compact manifold on
which $G$ acts smoothly. For any $\alpha \in H^3_G (M, \zz)$,
the equivariant Chern character composing with the isomorphism
induced by  \eqref{eq:tau}
leads to a natural  isomorphism

\begin{equation}
 K\upcom_{G, \alpha}(M)\otimes_{R(G)}R^\infty (G)
\iso  H_{G, delocalized, \alpha}\upcom (M)
\end{equation}
\end{conjecture}




\end{document}